\documentclass[smallextended,referee,envcountsect,]{svjour3}
\smartqed
\usepackage{graphicx}

\usepackage{amsmath}
\usepackage{amsfonts}
\usepackage{amssymb}
\usepackage{graphicx}
\usepackage{color}
\usepackage{url,algorithm}
\usepackage{enumitem} 
\usepackage{multirow}
\usepackage{tabularx}
\usepackage[title]{appendix}%
\usepackage{subfigure}

\newcommand{\rr}{{\mathbb R}}
\newcommand{\zz}{{\mathbb Z}}
\newcommand{\tr}{^{\mathsf T}}
\newcommand{\smallmat}[1]{\left[ \begin{smallmatrix}#1 \end{smallmatrix} \right]}

\newcommand{\ba}[1]{\begin{array}{#1}}
\newcommand{\ea}{\end{array}}

\begin{document}

\title{Global and Preference-based Optimization with Mixed Variables using Piecewise Affine Surrogates}


\author{Mengjia Zhu$^1$          \and
        Alberto Bemporad$^1$  
}

\institute{Mengjia Zhu, Corresponding Author \at
             mengjia.zhu@imtlucca.it
           \and
              Alberto Bemporad \at
              alberto.bemporad@imtlucca.it\\
              \and 
            $^1$ IMT School for Advanced Studies Lucca, Piazza San Francesco 19, 55100, Lucca, Italy\\
}

\date{Received: date / Accepted: date}

\maketitle

\begin{abstract}
Optimization problems involving mixed variables (\emph{i.e.}, variables of numerical and categorical nature) can be challenging to solve, especially in the presence of mixed-variable constraints. Moreover, when the objective function is the result of a complicated simulation or experiment, it may be expensive-to-evaluate. This paper proposes a novel surrogate-based global optimization algorithm to solve linearly constrained mixed-variable problems up to medium size (around 100 variables after encoding). The proposed approach is based on constructing a piecewise affine surrogate of the objective function over feasible samples. We assume the objective function is black-box and expensive-to-evaluate, while the linear constraints are quantifiable, unrelaxable, a priori known, and are cheap to evaluate. We introduce two types of exploration functions to efficiently search the feasible domain via mixed-integer linear programming solvers. We also provide a preference-based version of the algorithm designed for situations where only pairwise comparisons between samples can be acquired, while the underlying objective function to minimize remains unquantified. The two algorithms are evaluated on several unconstrained and constrained mixed-variable benchmark problems. The results show that, within a small number of required experiments/simulations, the proposed algorithms can often achieve better or comparable results than other existing methods.
\end{abstract}




\section{Introduction}\label{sec:intro}
A large variety of decision problems in several application domains, such as model selection in machine learning~\cite{luo2016review}, engineering design~\cite{kim2020surrogate}, and protein design~\cite{yang2019machine}, require identifying a global optimum without an explicit closed-form expression correlating the optimization variables to form the objective function to optimize. 
Such a black-box input/output function can be expensive-to-evaluate, as they may represent the outcomes (i.e., function outputs) of costly experiments or computationally intensive simulations resulting from a given value of the decision variables (i.e., function inputs).
Therefore, the number of function evaluations should be reduced as much as possible. 
Also, these functions often operate over mixed-variable domains, meaning the variables can be of different types, such as continuous, integer, or categorical, adding complexity to the optimization process. 
Additionally, physical problems frequently include constraints of a mixed-integer nature (e.g., constraints formed with logical conditions, which often involve continuous and binary variables). In some cases, evaluating infeasible instances of the optimization variables
may be impossible---for instance, when the corresponding function evaluation requires running a simulation or an experiment that can not be conducted or poses safety risks. As a result, it is preferable to efficiently exploit the known admissible set of the problem to encourage feasible sampling.

Surrogate-based optimization techniques have been studied extensively to target black-box optimization problems with expensive-to-evaluate objective functions~\cite{jones1998efficient,audet2017derivative}. For example, Bayesian Optimization (BO)~\cite{movckus1975bayesian,brochu2010tutorial} has been widely used in hyperparameter tuning in machine learning~\cite{wu2019hyperparameter} and adaptive experimental design~\cite{greenhill2020bayesian}. 
One of the authors recently developed a deterministic algorithm, GLIS~\cite{Bem20}, which has been applied to controller tuning~\cite{forgione2020efficient}, demonstrating 
comparable performance to BO but with lower computational cost. 
Although most of the literature has focused only on real-valued optimization variables, a few approaches have been adopted to handle integer and categorical variables~\cite{audet2023general}. Here, we distinguish integer variables as those representing ordinal relationships and categorical variables as those representing non-ordinal relationships. 
Integer variables are most commonly considered as continuous variables during the solution process and rounded to the nearest integer during post-analysis (\emph{e.g.}, MISO~\cite{muller2016miso}, RBFopt~\cite{costa2018rbfopt}). On the other hand, categorical variables are often one-hot encoded and treated as continuous variables within the range $[0,1]$ when constructing the surrogate model during optimization. After the optimization step, these variables are rounded and decoded back to their original categorical form for testing (\emph{e.g.}, One-hot BO~\cite{gpyopt2016}, MINOAN~\cite{kim2020surrogate}). See also~\cite{mistry2021mixed,holmstrom2008adaptive,gpyopt2016} for algorithms that have applied similar approaches to handle integer and categorical variables. 
Ploskas and Sahinidis~\cite{ploskas2022review} comprehensively analyzed and compared different algorithms and their associated software packages, targeting bound-constrained mixed-integer derivative-free optimization problems. In their review, the authors observed that MISO~\cite{muller2016miso} demonstrates superior performance when dealing with large (51 - 500 variables) and binary problems. On the other hand, NOMAD~\cite{abramson2011nomad,audet2021nomad} emerged as the top performer for mixed-integer, discrete (non-binary), small, and medium-sized (up to 50 variables) problems. 

Most of the surrogate-based methods assume all the inputs as continuous and ordinal~\cite{kim2020surrogate,muller2016miso,ploskas2022review}. On the other hand, different classes for the categorical variables often represent different choices rather than ordinal relations. Therefore, if one attempts to fit the latent function using a unified surrogate, in which categorical variables are one-hot encoded and treated as continuous vectors with entries in $[0,1]$,
sharp transitions might be observed in the constructed surrogate, leading to poor fitting qualities.
Alternatively, one can fit different surrogate models to each categorical class~\cite{swiler2014surrogate,gopakumar2018algorithmic,nguyen2019bayesian}. For example, EXP3BO~\cite{gopakumar2018algorithmic} constructs a surrogate model with Gaussian Process (GP) for each chosen class of the categorical variable.  However, as the number of categories and classes within each category increase, the size of the problem quickly blows up. To alleviate this issue, Ru et al.~\cite{ru2020bayesian} propose an approach that makes efficient use of the information in the acquired data by combining the strengths of multi-armed bandits and BO with GP models. This method, named as CoCaBO, has been shown to effectively solve bound-constrained problems with multiple categorical variables and multiple possible choices. 
Certain algorithms can inherently handle categorical variables due to the nature of their models. For example, the tree-structured models used in TPE~\cite{bergstra2011algorithms} and the random forests employed in SMAC~\cite{hutter2011sequential} naturally process categorical data.

In addition to mixed variables, real-life optimization problems frequently contain constraints. In this case, if the integer and the one-hot encoded categorical variables are relaxed as continuous variables while optimizing (\emph{i.e.}, the integrality of the variables is neglected) the constraints may not be satisfied after post-analysis, especially when equality constraints are present~\cite{kim2020surrogate}. In~\cite{kim2020surrogate}, to maintain the integrality of the variables, the authors use one-hot encoding to convert integer and one-hot encoded categorical variables to auxiliary variables. However, infeasibility with respect to constraints is still allowed during the solution process in~\cite{kim2020surrogate}.
In~\cite{papalexopoulos2022constrained}, piecewise-linear neural networks are employed as surrogate models to address constrained discrete black-box optimization problems, where mixed-integer linear programming (MILP) is used to optimize the acquisition function. However, the no-good constraints used in~\cite{papalexopoulos2022constrained} to tackle discrete-variable-only problems cannot be trivially transferred to the mixed-variable domain; hence, this approach cannot be directly applied to domains with mixed variables.

\subsection{Contribution}
In this work, we aim to solve medium sized mixed-variable nonlinear optimization problems (up to 100 variables after encoding) subject to mixed-integer linear equality and/or inequality constraints (up to 100 constraints).
Specifically, the optimization variables can be continuous, integer, and categorical, and the constraints are quantifiable unrelaxable a priori known (QUAK) based on the taxonomy in~\cite{le2024taxonomy}.
We
propose an algorithm that uses piecewise affine (PWA) functions as the surrogate models. PWA functions can effectively handle the discontinuities introduced by sharp transitions among different classes of categorical variables. Furthermore, they can be directly reformulated into MILPs, allowing us to leverage efficient, off-the-shelf MILP solvers to optimize the acquisition function.
To balance the exploitation and exploration during acquisition, we incorporate two types of exploration functions---distance-based and frequency-based---in the acquisition function.  
Additionally, we propose incorporating the exploration function as part of the initial sampling strategy to obtain well-scattered initial samples, especially when a large number of linear equality and/or inequality constraints are present. This is crucial because the initial samples play an essential role in surrogate fitting, particularly when the function-evaluation budget is limited.


We name the proposed algorithm as \textbf{PWAS}, short for \textbf{P}iece\textbf{w}ise \textbf{A}ffine \textbf{S}urrogate-based optimization. We show the efficiency and effectiveness of PWAS by comparing its performance with other existing solvers on a set of benchmark problems. We also present an extension of PWAS to address problems where function evaluations are unavailable, such as those involving multiple objectives with unclear relative weights to form a single objective function or cases where only qualitative assessments are available. This approach assumes that a decision-maker can express \emph{preferences} between two candidate solution vectors.
Such preference information is used to shape the PWA surrogate through the proposed algorithm, named \textbf{PWASp} (short for \textbf{PWAS} based on \textbf{p}references).
Python implementations of PWAS and PWASp are available on the GitHub repository (\url{https://GitHub.com/mjzhu-p/PWAS}).

The rest of the paper is organized as follows. The description of the target problem is formulated in Section~\ref{sec:prob_decp}. The proposed surrogate-based optimization algorithms are discussed in Sections~\ref{sec:prop_med} and~\ref{sec:pref_learning}. Section~\ref{sec:benchmark} reports the numerical benchmarks demonstrating the effectiveness of the proposed method. Lastly, conclusions and directions for future research are discussed in Section~\ref{sec:conclusion}.

\section{Problem Formulation}\label{sec:prob_decp}
We consider a decision problem with $n_c$ real variables
grouped in vector ${x\in\rr^{n_c}}$, $n_{\rm int}$ integer variables
grouped in vector $y\in\zz^{n_{\rm int}}$, and $n_d$ categorical variables grouped in list $Z = [Z^1,\ldots,Z^{n_d}]$. Each categorical variable $Z^i$ can
take values within its corresponding $n_i$ classes, for $i=1,\ldots,n_d$. 
Let us assume that each categorical variable $Z^i$ is one-hot binary encoded
into the subvector $[z_{1+d^{i-1}}\ \ldots\ z_{d^i}]\tr\in\{0,1\}^{n_i}$
for each $i=1,\ldots,n_d$, where $d^0=0$ and $d^i=\sum_{j=1}^{i} n_j$, with $z\in\{0,1\}^{d^{n_d}}$ being the complete vector of binary variables
after the encoding. Here, $z\in\Omega_z=\{z\in\{0,1\}^{d^{n_d}}:\
\sum_{j=1}^{n_i} z_{j+d^{i-1}}=1,\ \forall i=1,\ldots,n_d\}$.
Let $X=[x\tr\ y\tr\ z\tr]\tr$ denote the overall optimization vector. 
We assume that the vectors $x$ and $y$ of interest 
are bounded (\emph{i.e.}, $\ell_x \leq x \leq u_x$ and $\ell_y \leq y \leq u_y$)
and denote 
the domain of $X$ by $\Omega=[\ell_x,u_x]\times([\ell_y,u_y]\cap \zz)\times\Omega_z$. 
We denote $f:\Omega\mapsto\rr$ as the objective function to minimize, assuming it is noiseless and expensive-to-evaluate. 

The black-box mixed-variable optimization problem we want to solve can be stated as follows, where we aim to 
\begin{subequations}
\label{eq:opt_prob_1}
\begin{eqnarray}
     {\rm find}\ \ X^* \in&& \arg \min_{X \in \Omega} f(X)\\
  \text{s.t.}  &&\quad A_{\rm eq}x + B_{\rm eq}y + C_{\rm eq}z= b_{\rm eq}\label{eq:opt_prob_1-b}\\
  			&&\quad A_{\rm ineq}x + B_{\rm ineq}y + C_{\rm ineq}z \leq b_{\rm ineq},
\label{eq:opt_prob_1-c}
\end{eqnarray}
\end{subequations}
where $A_{\rm eq}$, $B_{\rm eq}$, $C_{\rm eq}$, $A_{\rm ineq}$, $B_{\rm ineq}$, and $C_{\rm ineq}$ are matrices of suitable dimensions that, together with 
the right-hand-side column vectors $b_{\rm eq}$ and $b_{\rm ineq}$, define possible
linear equality and inequality constraints on $x$, $y$, and $z$. For example, if $x\in\rr$ and $Z = [Z^1]$ with $Z^1\in\{\textsf{red},\textsf{blue},\textsf{yellow}\}$, the logical constraint $[Z^1=\textsf{red}]\rightarrow [x\leq 0]$ can be modeled as
$x\leq u_x (1-z_1)$. For modeling more general types of mixed linear/logical constraints (possibly involving the addition of auxiliary real and binary variables), the reader is referred to, e.g.,~\cite{Wil13,HO99,TB04a}. Note that, different from function $f$, which is assumed to be black-box and expensive-to-evaluate, we assume the mixed-integer linear constraints on $X$ in~\eqref{eq:opt_prob_1} are QUAK~\cite{le2024taxonomy} and are cheap to evaluate. 

\section{Solution Method}\label{sec:prop_med}
\begin{figure}[t]
    \centering
    \includegraphics[width =\textwidth]{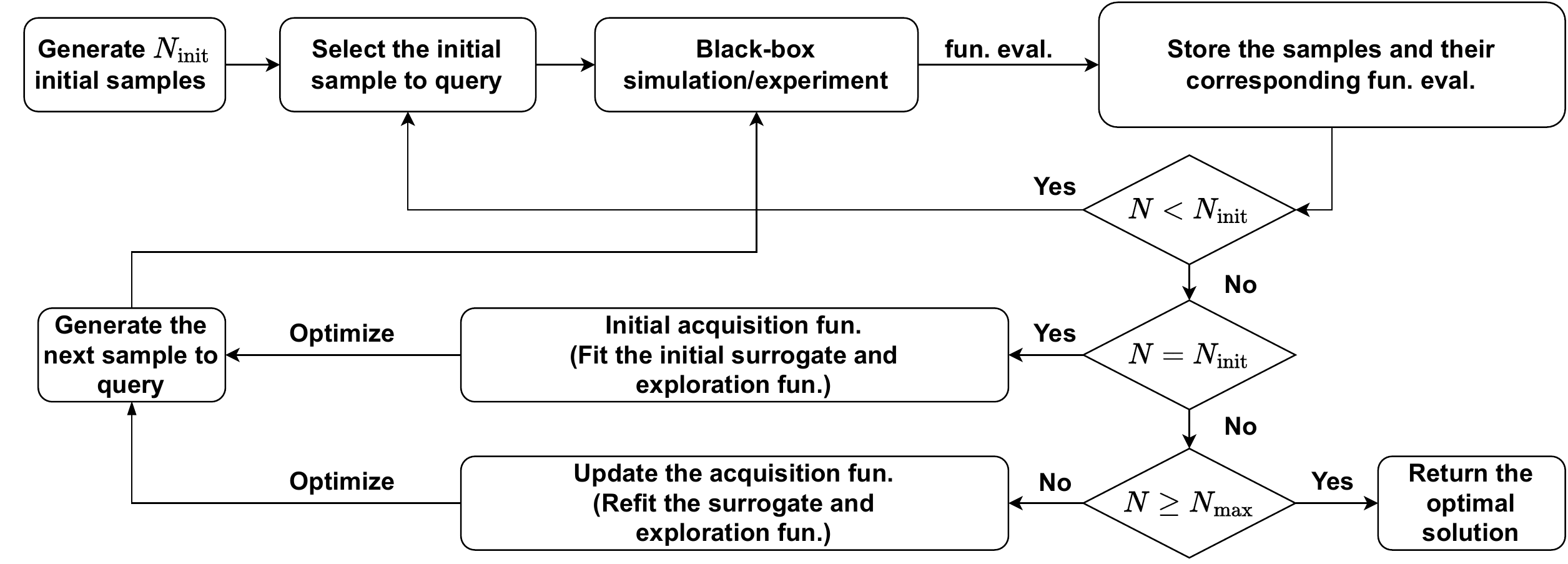}
    \caption{General procedures for surrogate-based optimization methods.}
    \label{fig:surrogate_based_opt_general_flow}
\end{figure}
We follow the general surrogate-based optimization procedure, schematically depicted in Figure~\ref{fig:surrogate_based_opt_general_flow} (see, \emph{e.g.},~\cite{Bem20}), to solve the optimization problem~\eqref{eq:opt_prob_1}. The approach consists of an initial (passive) sampling and an active learning stage, in which a surrogate model of the objective function $f$ is repeatedly learned. During the initial phase, $N_{\rm init}$ feasible samples $X_k$ are generated, and the corresponding values $f(X_k)$ are collected, for $k=1,\ldots,N_{\rm init}$. During active learning, the surrogate model of $f$ 
is estimated from a finite number of function evaluations $f(X_k)$, for $k=1,\ldots,N$, where $N$ is increased sequentially between the initial value $N_{\rm init}$ and
the maximum budget $N_{\rm max}$ of queries available. Each new sample $X_{k+1}$
is determined by minimizing an \emph{acquisition function}, which combines the surrogate with an \emph{exploration function}, to reach an exploitation/exploration tradeoff. 
The exploration function is constructed based on the existing samples $\{X_k\}$, for $k=1,\ldots,N$, to ensure that the feasible domain is sufficiently explored. The procedure aims to effectively reduce the objective function value of the posed problem within a small number of function evaluations.

In this paper, we propose fitting a \emph{piecewise affine} (PWA) surrogate of the latent objective function $f$. PWA surrogate has two main benefits: ($i$) allow discontinuities introduced by sharp transitions induced by taking values in different classes of the categorical variables. 
In this case, instead of using one surrogate model for each categorical class as in~\cite{swiler2014surrogate,gopakumar2018algorithmic,nguyen2019bayesian}, it is possible to adaptively update the number of partitions allowed in the PWA function by analyzing the clusters of the queried samples. For example, one can initiate the surrogate fitting procedure by setting a maximum allowed number of partitions and then discard some partitions if the number of queried samples within these partitions is smaller than some fixed minimum values (cf.~\cite{BemPARC21}); ($ii$) PWA surrogates have a direct mixed-integer linear reformulation and, therefore, can be minimized by efficient MILP solvers (\emph{e.g.}, Gurobi~\cite{gurobi} and GLPK~\cite{makhorin2020gnu}). 
Also, we can explicitly reformulate and include linear equality and inequality constraints involving integer and one-hot encoded categorical variables in the standard MILP form to maintain their integrality during the solution process, enabling the possibility to make feasible queries during the acquisition step.
Since the acquisition function includes both the \emph{surrogate} and the \emph{exploration function}, we will also define a suitable PWA \emph{exploration function} that admits a MILP representation. The resulting approach, which we call PWAS, is summarized in Algorithm~\ref{algo:pwas}, whose steps will be described in detail in the next sections.

\begin{algorithm}[ht!]
	\caption{PWAS: Global Optimization Using Piecewise Affine Surrogates}
	\label{algo:pwas}
	~~\textbf{Input}: Lower and upper bounds $\ell_x,u_x,\ell_y,u_y$; linear constraint matrices $A_{\rm eq}$, $B_{\rm eq}$, $C_{\rm eq}$, $A_{\rm ineq}$, $B_{\rm ineq}$ and $C_{\rm ineq}$; linear constraint right-hand-side vectors $b_{\rm eq}$ and $b_{\rm ineq}$; number $n_d$ of categorical variables and $n_i$ of possible categories, for $i=1,\ldots,n_d$; initial number $K$ of polyhedral partitions;
	number $N_{\rm init}\geq 2$ of initial samples; number $N_{\rm max}\geq N_{\rm init}$ of maximum function evaluations; $\delta_1 \geq 0$, $\delta_2 \geq 0$ and $\delta_3 \geq 0$ if solve~\eqref{eq:acq_solve} in one step or $\delta \geq 0$ if solve~\eqref{eq:acq_solve} in multiple steps; solving strategy for~\eqref{eq:acq_solve}: \{``one-step'' or ``multi-steps''\}.
	\vspace*{.1cm}\hrule\vspace*{.1cm}
	\begin{enumerate}
		\item Pre-process the optimization variables as described in Section~\ref{subsec:pre-processing};
		\item  $N\leftarrow 1$, $N^*_{\rm curr}\leftarrow 1$, $f^*_{\rm curr}\leftarrow +\infty$ ;
		\item \label{algo:initial} Generate $N_{\rm init}$ random scaled and encoded samples $\bar X=\{\bar X_1,\ldots, \bar X_{N_{\rm init}}\}$ using one of the initial sampling methods reported in Section~\ref{subsec:initial_sampl} based on the problem setup;

			\item \textbf{While} $N \leq N_{\rm max}$ \textbf{do}
			\begin{enumerate}
				\item Scale back and decode $\bar X_N$ to $X_N$, \emph{i.e.}, $X_N = S(\bar X_N)$, and query $f_N=f(X_N)$;
				\item \textbf{If} $f_N < f^*_{\rm curr}$ \textbf{then update}  $N^*_{\rm curr}\leftarrow N$, $f^*_{\rm curr}\leftarrow f_N$;
				\item \textbf{If} $N \geq N_{\rm init}$ \textbf{then}
				\begin{enumerate}
					\item Update and fit the PWA separation function $\phi$ and PWA surrogate function $\hat f$ as described in Section~\ref{subsec:PWA_surrogate};
					\item Define the acquisition function $a$ as in~\eqref{eq:acq_solve};
					\item Solve the global optimization problem~\eqref{eq:acq_solve} and get $\bar X_{N+1}$ either in one-step or multi-steps;
				\end{enumerate}
				\item $N\leftarrow N + 1$;
			\end{enumerate}
		\item \textbf{End}.
	\end{enumerate}
	\vspace*{.1cm}\hrule\vspace*{.1cm}
	~~\textbf{Output}: Best decision vector $X^*=X_{N^*_{\rm curr}}$ found.
\end{algorithm}

\subsection{Change of Variables: Scaling and Encoding}\label{subsec:pre-processing}
Before attempting solving problem~\eqref{eq:opt_prob_1}, we first rescale every continuous variable $x_i$ into a new variable $\bar x_i\in[-1, 1]$ such that
\[
	x_i =  \frac{u_x^i - \ell_x^i}{2} \bar x_i + \frac{ u_x^i +\ell_x^i}{2},\ \forall i = 1,\dots, n_c.
\]
Accordingly, the constraint matrices $A_{\rm eq}$ and $A_{\rm ineq}$ are rescaled to
\[
	\bar A_{\rm eq} = A_{\rm eq} {\rm diag}\left( \frac{ u_x-\ell_{x}}{2} \right)\ {\rm and}\  \bar A_{\rm ineq} = A_{\rm ineq} {\rm diag}\left( \frac{ u_{x}-\ell_{x}}{2} \right),
\]
with the right-hand-side vectors updated as follows:
\[
	\bar b_{\rm eq} = b_{\rm eq} - A_{\rm eq}\left( \frac{ u_x+\ell_{x}}{2} \right)\ {\rm and}\  \bar b_{\rm ineq} = b_{\rm ineq} - A_{\rm ineq}\left( \frac{ u_x+\ell_{x}}{2} \right).
\]
The intervals $[-1,1]$ for the continuous variables are possibly further tightened by taking the updated inequality constraints~\eqref{eq:opt_prob_1-c} (if they exist) into account (cf.~\cite{Bem20}), \emph{i.e.}, for $i=1,\ldots,n_c$, we set 
\[
\begin{array}{rl}
	\bar\ell_x^i &= \underset{\bar x, y, z}{\min} \  e_{i} \tr [\bar x \tr\  y \tr\  z \tr]\tr\\ [1em]
	\text{s.t.} &\quad \bar A_{\rm ineq}\bar x + B_{\rm ineq}y + C_{\rm ineq}z \leq \bar b_{\rm ineq}\\
  & \quad \bar x\in[-1\  1]^{n_c},\   y \in [\ell_y,u_y]\cap \zz,\  z \in \Omega_z,
\end{array}
\]
and, similarly,
\[
\begin{array}{rl}
	\bar u_x^i &= \underset{\bar x, y, z}{\max} \  e_{i} \tr [\bar x \tr\  y \tr\  z \tr]\tr\\ [1em]
	\text{s.t.} &\quad \bar A_{\rm ineq}\bar x + B_{\rm ineq}y + C_{\rm ineq}z \leq \bar b_{\rm ineq}\\
  & \quad \bar x\in[-1\  1]^{n_c},\   y \in [\ell_y,u_y]\cap \zz,\  z \in \Omega_z.
  \end{array}
\]
Here, $e_i$ denotes the $i$th column of the identity matrix of the same dimension
as vector $X$. We denote the resulting domain of the scaled continuous variables by $\Omega_x=[\bar \ell_x^1,\bar u_x^1]\times\ldots\times[\bar \ell_x^{n_c},\bar u_x^{n_c}]$.

Let us assume that only a finite number $N_{\rm max}$ of queries can be made,
which depends on the nature of function $f$ (i.e., how expensive it is to evaluate) and the time available to solve the optimization problem. 
Moreover, we treat integer variables $y$ differently
depending on the relation between $N_{\rm max}$ 
and the number $\prod_{i = 1}^{n_{\rm int}} n_i^{\rm int}$
of possible combinations of integer variables, where $n_i^{\rm int} = \lfloor u^i_y\rfloor-\lceil \ell^i_y \rceil +1$ is the cardinality of the set $[l_y^i,u_y^i]\cap\zz$, \emph{i.e.}, the number of integer values that variable $y_i$ can take.
To be described in Section~\ref{sec:scenario1} and~\ref{sec:scenario2} below, we will treat the vector $y$ of integer variables in two distinct ways to make the exploration of the search space possibly more efficient. Specifically, $y$ is treated as categorical when solving~\eqref{eq:opt_prob_1} in case $\prod_{i = 1}^{n_{\rm int}} n_i^{\rm int}<N_{\rm max}$, \emph{i.e.}, when it may be possible to exhaustively list out all the potential combinations of the integer variables within $N_{\rm max}$ queries if no continuous or categorical variables are present; vice versa, we will maintain the optimization variables $y_i$ integer. We note that this is a general heuristic we applied, which was empirically observed to be more efficient when handling integer variables. This heuristic is motivated by the acquisition strategies, which we will elaborate more in Section~\ref{subsec:acq_fun}.

\subsubsection{Treating Integer Variables as Categorical}
\label{sec:scenario1}
The first scenario occurs when the number of possible combinations of integer variables $\prod_{i = 1}^{n_{\rm int}} n_i^{\rm int}<N_{\rm max}$.
In this case, we treat all integer variables $y_i$ as categorical (similarly to vector $z$) and one-hot encode them into further $d^{n_{\rm int}}$ binary variables $\bar y_j\in\{0,1\}$, for $j=1,\ldots,d^{n_{\rm int}}$, where $d^{n_{\rm int}}=\sum_{i=1}^{n_{\rm int}}n_i^{\rm int}$. We also define $\Omega_y=\{\bar y\in\{0,1\}^{d^{n_{\rm int}}}:\
\sum_{j=1}^{n_i^{\rm int}} \bar y_{j+d_y^{i-1}}=1,\ \forall i=1,\ldots,n_{\rm int}\}$,
where $d_y^0=0$ and $d_y^i=\sum_{j=1}^{i}n_j^{\rm int}$ for $i=1,\ldots,n_{\rm int}$,
and set $\bar y\in\Omega_y$. The constraint matrix $B_{\rm eq}$ ($B_{\rm ineq}$) is modified accordingly into a new matrix $\bar B_{\rm eq}$ ($\bar B_{\rm ineq}$) by replacing each scalar entry $B_{\rm eq}^{ij}$ ($B_{\rm ineq}^{ij}$) with the
row vector obtained by multiplying the entry by the vector of integers
that variable $y_j$ can take, \emph{i.e.},
\[
	\begin{array}{rcll}
	B_{\rm eq}^{ij} &\leftarrow &B_{\rm eq}^{ij} \left[\lceil \ell^j_y \rceil\ \ldots\ \lfloor u^j_y\rfloor \right]\in\rr^{1\times n_j^{\rm int}}, &\forall j = 1, \dots, n_{\rm int}, \\
	B_{\rm ineq}^{ij} &\leftarrow &B_{\rm ineq}^{ij} \left[\lceil \ell^j_y \rceil\ \ldots\ \lfloor u^j_y\rfloor  \right]\in\rr^{1\times n_j^{\rm int}}, &\forall j = 1, \dots, n_{\rm int}.
	\end{array}
\]
The new optimization vector becomes $\bar X=[\bar x\tr\ \bar y\tr\ z\tr]\tr\in\bar\Omega$, where $\bar\Omega=\Omega_x\times\Omega_y\times\Omega_z$, consisting of $n=n_c+d^{n_{\rm int}}+d^{n_d}$ variables. As evaluating the objective function in~\eqref{eq:opt_prob_1} requires the original values in $X$, we denote 
by $S:\bar\Omega\mapsto\Omega$ the inverse scaling/encoding mapping of $\bar X$, \emph{i.e.}, $X=S(\bar X)$. According to such a change of variables,
problem~\eqref{eq:opt_prob_1} is now translated to 
\begin{subequations}
\label{eq:opt_prob_2}%
\begin{equation}
	\begin{split}
    \text{find}\ \bar X^* \in\ & \arg \min_{\bar X \in \bar \Omega} f(S(\bar X))\\
  \text{s.t.}\  
  &\bar A_{\rm eq} \bar x + \bar B_{\rm eq}\bar y + C_{\rm eq}z= \bar b_{\rm eq}\\
&\bar A_{\rm ineq} \bar x + \bar B_{\rm ineq} \bar y + C_{\rm ineq}z \leq \bar b_{\rm ineq}.
\end{split}
\label{eq:opt_prob_2-a}
\end{equation}
In the sequel, $D\subseteq\bar \Omega$ will denote the set of admissible 
vectors $\bar X$ satisfying the constraints in~\eqref{eq:opt_prob_2-a}.

\subsubsection{Scaling Integer Variables}
\label{sec:scenario2}
In the second scenario where $\prod_{i = 1}^{n_{\rm int}} n_i^{\rm int}\geq N_{\rm max}$, the integer variables are rescaled and treated as numeric variables
$\bar y_i\in[-1,1]$, $i=1,\ldots,n_{\rm int}$. In this case, 
we also keep the original $n_{\rm int}$ integer variables $y_i\in\zz$ in the model 
for the sole purpose of enforcing integrality constraints, as we 
link them with $\bar y_i$ by the scaling factors 
\[
	y_i =  \frac{ u_y^{i} - \ell_y^{i} }{2} \bar y_i + \frac{ u_y^{i}  +\ell_y^{i} }{2}.  \\
\]
Similar to the continuous variables, we can also further shrink the bounds on $\bar y_i$ by considering the updated inequality constraints~\eqref{eq:opt_prob_1-c} (if present)
\[
\begin{array}{rl}
	\bar\ell_y^i &= \underset{\bar x, \bar y, z}{\min} \  e_{n_c+i} \tr [\bar x \tr\  \bar y \tr\  z \tr]\tr\\ [1em]
	\text{s.t.} &\quad \bar A_{\rm ineq}\bar x + B_{\rm ineq}y + C_{\rm ineq}z \leq \bar b_{\rm ineq}\\
  & \quad \bar x\in \Omega_x,\   \bar y \in [-1\ 1]^{n_{\rm int}},\  z \in \Omega_z,
\end{array}
\]
\[
\begin{array}{rl}
	\bar u_y^i &= \underset{\bar x, \bar y, z}{\max} \  e_{n_c+i} \tr [\bar x \tr\  \bar y \tr\  z \tr]\tr\\ [1em]
	\text{s.t.} &\quad \bar A_{\rm ineq}\bar x + B_{\rm ineq}y + C_{\rm ineq}z \leq \bar b_{\rm ineq}\\
  & \quad \bar x\in \Omega_x,\   y \in[-1\ 1]^{n_{\rm int}},\  z \in \Omega_z.
  \end{array}
\]
We denote the domain of $\bar y$ after tightening as $\Omega_y=[\bar \ell_y^1,\bar u_y^1]\times\ldots\times[\bar \ell_y^{n_{\rm int}},\bar u_y^{n_{\rm int}}]$.
Accordingly, problem~\eqref{eq:opt_prob_1} is translated to
\begin{equation}
	\begin{split}
    \text{find}\ \begin{bmatrix}\bar X^*\\y^*\end{bmatrix} \in\ & \arg \min_{\bar X \in \bar \Omega, y\in[\ell_y,u_y]\cap \zz} f(S(\bar X))\\
  \text{s.t.}\  
  &\bar A_{\rm eq} \bar x + B_{\rm eq} y + C_{\rm eq}z= \bar b_{\rm eq}\\
&\bar A_{\rm ineq} \bar x + B_{\rm ineq} y + C_{\rm ineq}z \leq \bar b_{\rm ineq},
\end{split}
\label{eq:opt_prob_2-b}
\end{equation}%
\end{subequations}
where now $\bar X=[\bar x\tr\ \bar y\tr\ z\tr]\tr\in\bar\Omega$ and consists of $n=n_c+n_{\rm int}+d^{n_d}$ variables, with $\bar\Omega=\Omega_x\times\Omega_y\times\Omega_z$. And $S:\bar\Omega\mapsto\Omega$ is the new inverse scaling
mapping. 
We will denote by $D\subseteq\bar \Omega$ the set of admissible 
vectors $\bar X$ such that the constraints in~\eqref{eq:opt_prob_2-b}
are satisfied for some vector $y\in[\ell_y,u_y]\cap \zz$. 


\subsection{Piecewise Affine Surrogate Function}\label{subsec:PWA_surrogate}
When fitting a surrogate of the objective function, we treat the modified vector $\bar X$ as a vector in $\rr^n$. We describe next how to construct a PWA surrogate function $\hat f:\rr^n\mapsto\rr$ such that $\hat f(\bar X)$ approximates $f(S(\bar X))$.

Consider $N$ samples $\bar X_1$, $\ldots$, $\bar X_N\in\rr^{n}$ and their corresponding function evaluations $f(S(\bar X_1))$, $\ldots$, $f(S(\bar X_N)) \in \rr$.
We want to define the PWA surrogate function $\hat f$ over a polyhedral partition of $\bar\Omega$ into $K$ regions. To this end, we consider the following convex \emph{PWA separation function} $\phi:\rr^n\mapsto\rr$ 
\begin{subequations}
\begin{equation}\label{eq:pwa_partition}
   		\phi(\bar X)=\omega_{j(\bar X)}\tr\bar X+\gamma_{j(\bar X)},\
\end{equation}
where $\omega_j\in\rr^n$ and $\gamma_j\in\rr$, for $j=1,\ldots,K$ are the parameters that need to be determined,
with
\begin{equation}
	j(\bar X)=\arg\max_{j=1,\ldots,K}\{\omega_j\tr\bar X+\gamma_j\}.
	\label{eq:j(X)}
\end{equation}
We define the \emph{PWA surrogate function} $\hat f$ as
\begin{equation} \label{eq:pwa_eqn}
       \hat f(\bar X)=a_{j(\bar X)}\tr\bar X+b_{j(\bar X)},
\end{equation}
\end{subequations}
where $a_j\in\rr^n$ and $b_j\in\rr$, for $j=1,\ldots,K$, are the surrogate parameters that need to be determined.
Note that $\hat f$ is possibly non-convex and discontinuous. 

We use the PARC algorithm recently proposed in~\cite{BemPARC21} by one of the authors to fit the PWA separation and surrogate functions to obtain the required coefficients $\omega_j, \gamma_j, a_j, \ {\rm and}\ b_j$, for $j = 1, \dots, K$. We stress that while the closed-form expression of $f\circ S$ as a function of $\bar X$ is generally unavailable and very expensive-to-evaluate for each given $\bar X$, evaluating its surrogate $\hat f$ is very cheap
and, as we will show in Section~\ref{sec:MILP-surrogate}, admits a simple mixed-integer
linear encoding with $K$ binary variables.

The number of partitions, $K$, is a hyper-parameter of the proposed
global optimization algorithm that must be selected by trading
off between having a more flexible surrogate function (large $K$)
and reducing computational demands (small $K$). The PARC algorithm includes an adaptive mechanism to update $K$ during surrogate fitting, specifically, if a partition contains fewer samples than a predefined minimum, it is discarded. The samples from these discarded partitions are then either reassigned to a neighboring partition or treated as outliers and excluded from surrogate fitting, depending on their function evaluations~\cite{BemPARC21}. This adaptability makes the PARC algorithm more robust and flexible, allowing it to reduce computational complexity without significantly affecting prediction accuracy. For a detailed analysis of the PARC algorithm, we refer readers to~\cite{BemPARC21}. For illustration purposes, we present two surrogate-fitting examples: one in a continuous domain and another in a mixed continuous and categorical domain. While these demonstrations use 2D problems for simplicity and better visualization, the PARC algorithm is capable of handling high-dimensional problems.

\begin{figure}[tbh!]
    \centering
    \includegraphics[width=0.55\textwidth]{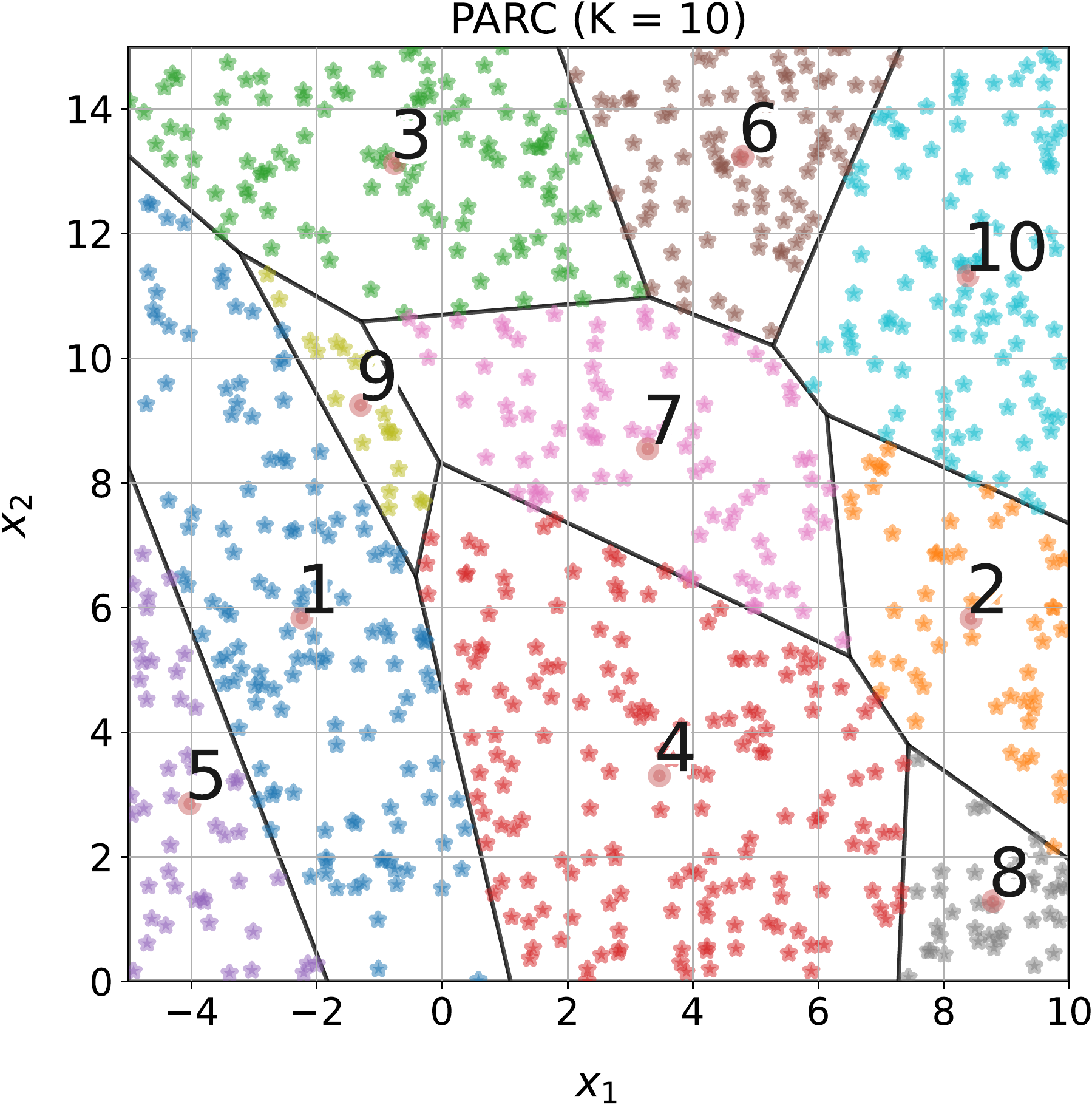}
    \caption{The polyhedral partition induced by~\eqref{eq:pwa_partition} with $K = 10$ initial partitions. Note that the final partition is also 10 (no partition is discarded). The dots in the figure are the training data for the PARC algorithm.Different colors indicate samples at different partitions. Brownish dots next to the partition numbers are the centroid of the training data within each partition.}
    \label{fig:parc_partition}
\end{figure}

\begin{figure}[tbh!]
    \centering
    \subfigure[]{\includegraphics[width=0.45\textwidth]{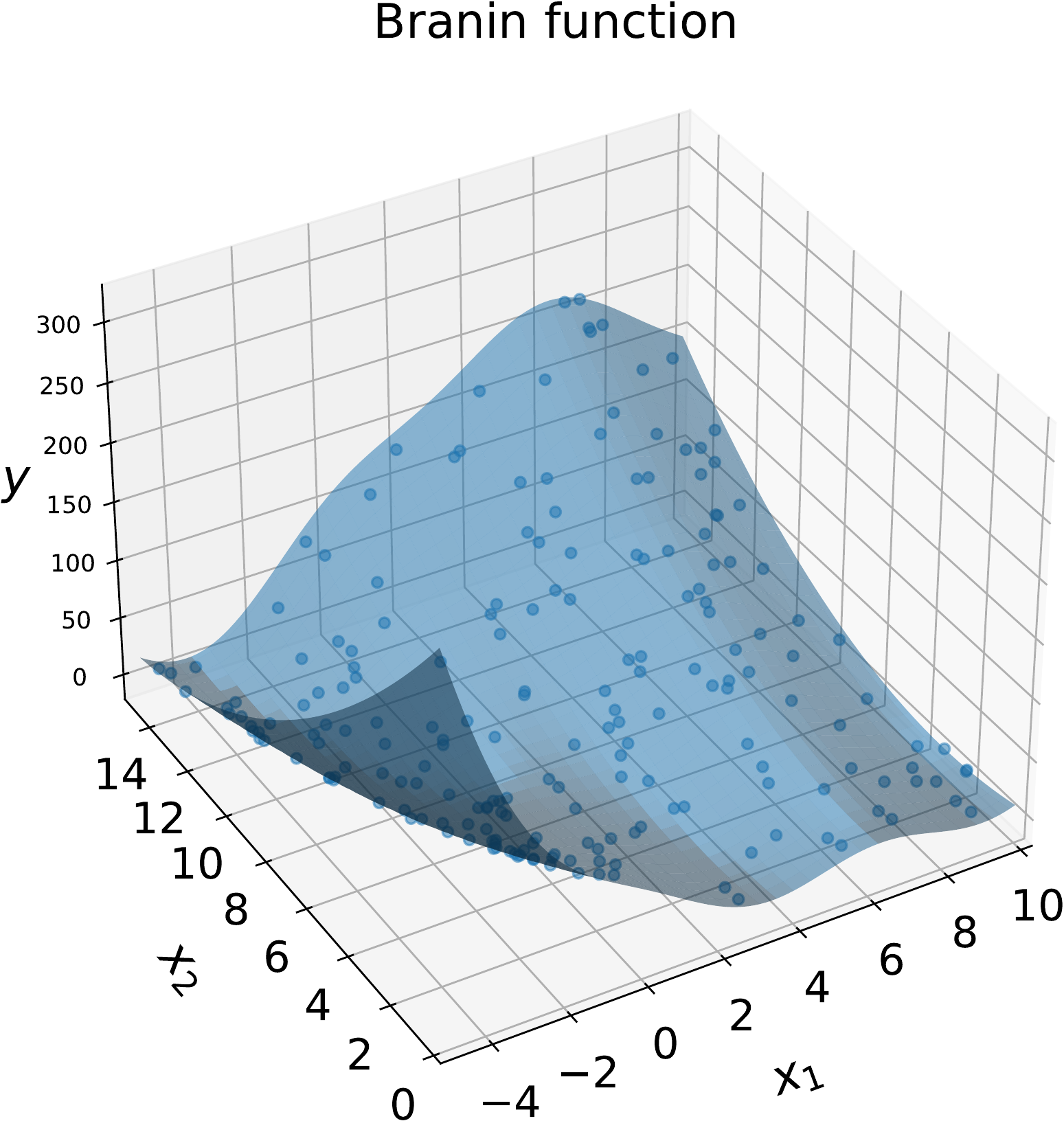}}\quad\quad 
    \subfigure[]{\includegraphics[width=0.45\textwidth]{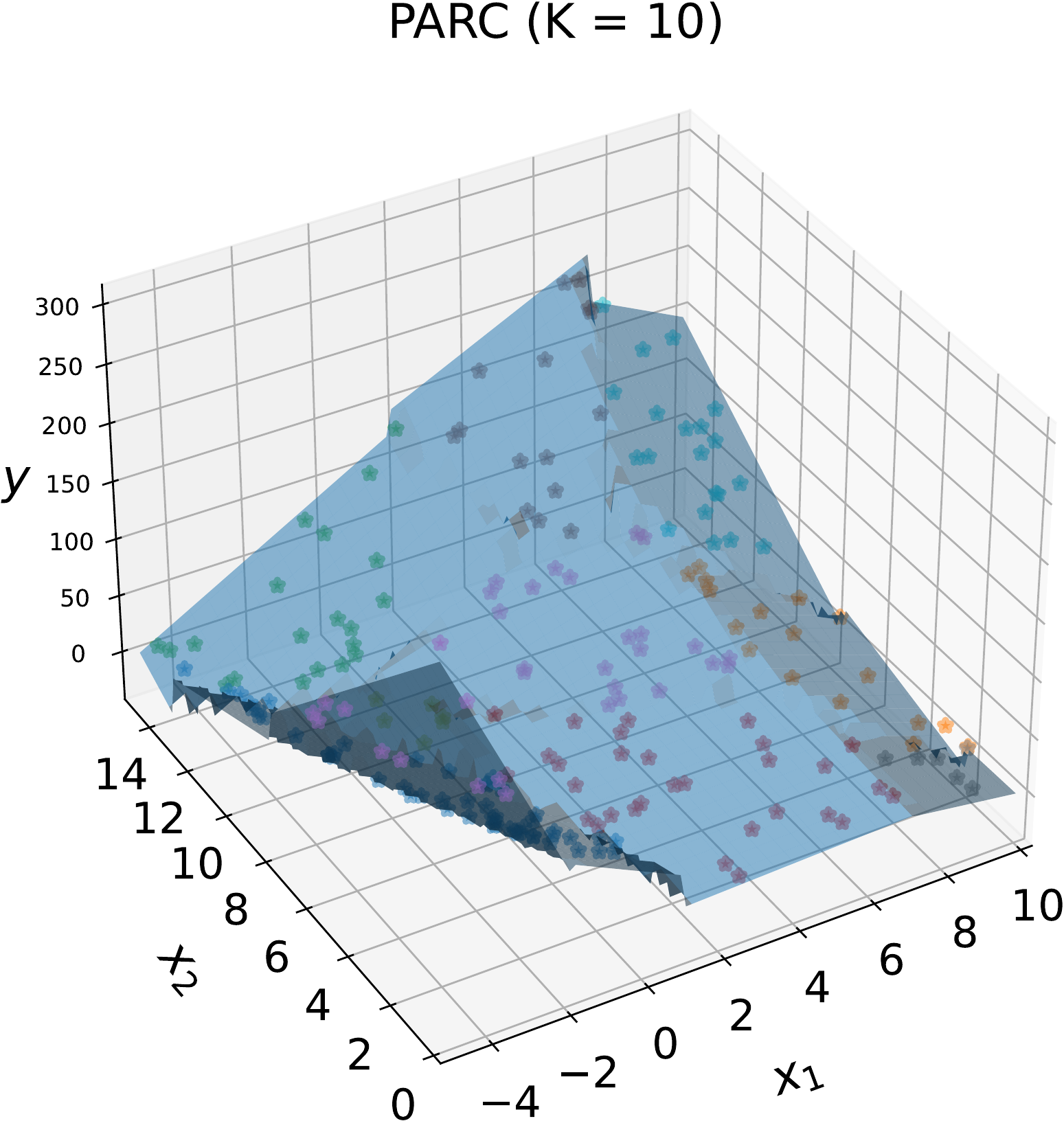}}
    \caption{(a) Branin function - analytical (b) Branin function fitted by PARC with $K = 10$. The dots in the figure are the test data (200 samples) for the PARC algorithm. Different colors in (b) indicate samples at different partitions.}
    \label{fig:pwe_fun_fitting}
\end{figure}

For the continuous function, we consider the Branin function~\cite{dixon1978global}:
\begin{equation}\label{eq:brainin_fun}
\begin{split}
    f(x_1,x_2) &= a(x_2 - bx_1^2+cx_1-r)^2 + s(1-t)\cos(x_1) +s \\
    a  &= 1,\  b= \frac{5.1}{4\pi^2}, \ c = \frac{5}{\pi},\  r = 6, \ s = 10,\  t = \frac{1}{8\pi}\\
    -5 & \leq x_1 \leq 10, \ 0 \leq x_2 \leq 15.
\end{split}
\end{equation}
A PWA surrogate of $f$ in~\eqref{eq:brainin_fun} is fitted by the PARC algorithm utilizing 800 randomly generated training samples $(x_{1k},x_{2k})$. The initial partition $K$ was set to 10. Figure~\ref{fig:parc_partition} show the final polyhedral partition induced by~\eqref{eq:pwa_partition}. And Figure~\ref{fig:pwe_fun_fitting} depicts the Branin function as represented both analytically and through the PWA surrogate fitted by PARC. We observe that the surrogate fitted by PARC 
captures the general shape of the nonlinear Branin function. 

\begin{figure}[tbh!]
    \centering
    \includegraphics[width=0.55\textwidth]{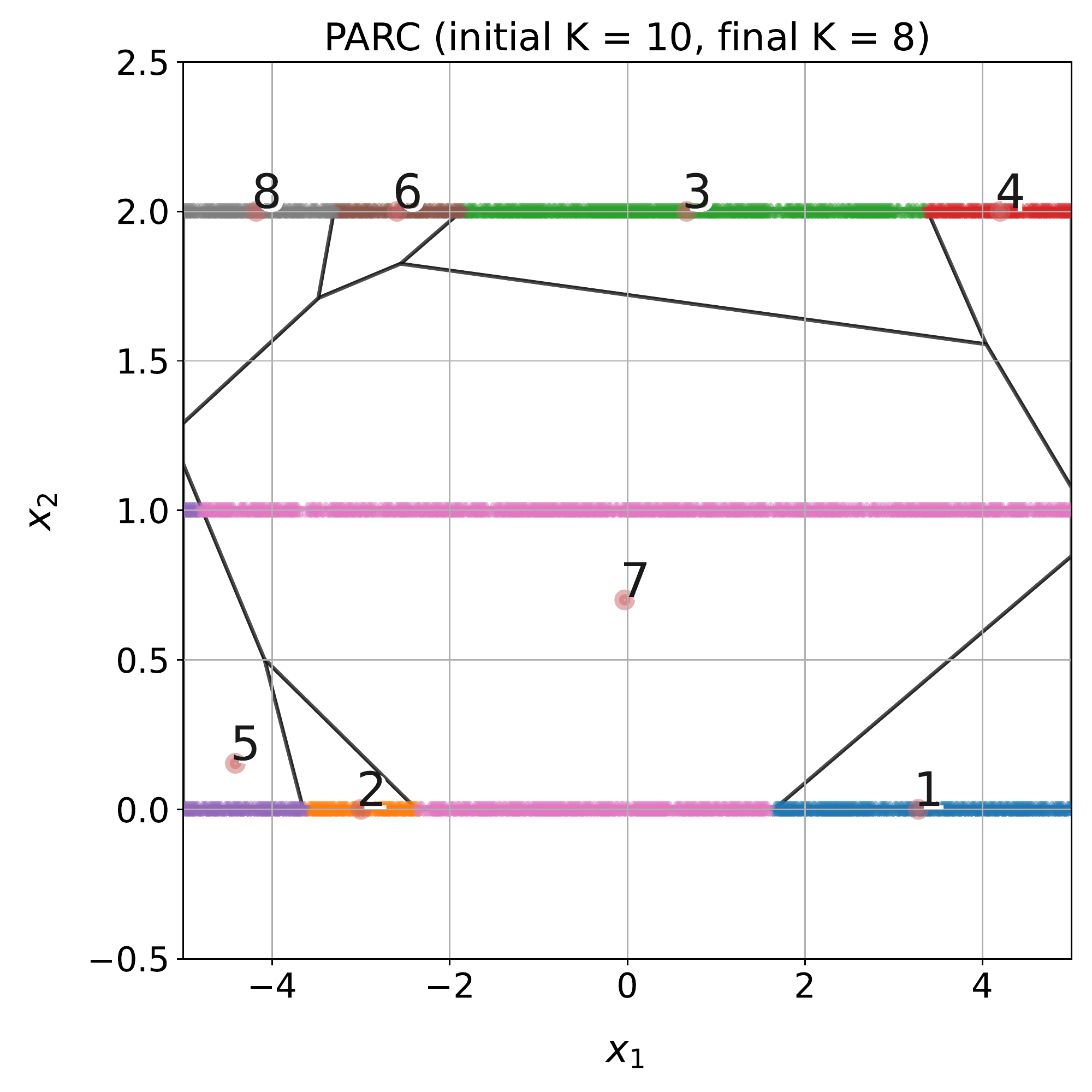}
    \caption{The polyhedral partition induced by~\eqref{eq:pwa_partition} with $K = 10$ initial partitions. Note that the final partition is 8 (2 partitions are discarded). The dots in the figure are the training data for the PARC algorithm.Different colors indicate samples at different partitions. Brownish dots next to the partition numbers are the centroid of the training data within each partition.}
    \label{fig:parc_partition_2}
\end{figure}

\begin{figure}[tbh!]
    \centering
    \subfigure[]{\includegraphics[width=0.45\textwidth]{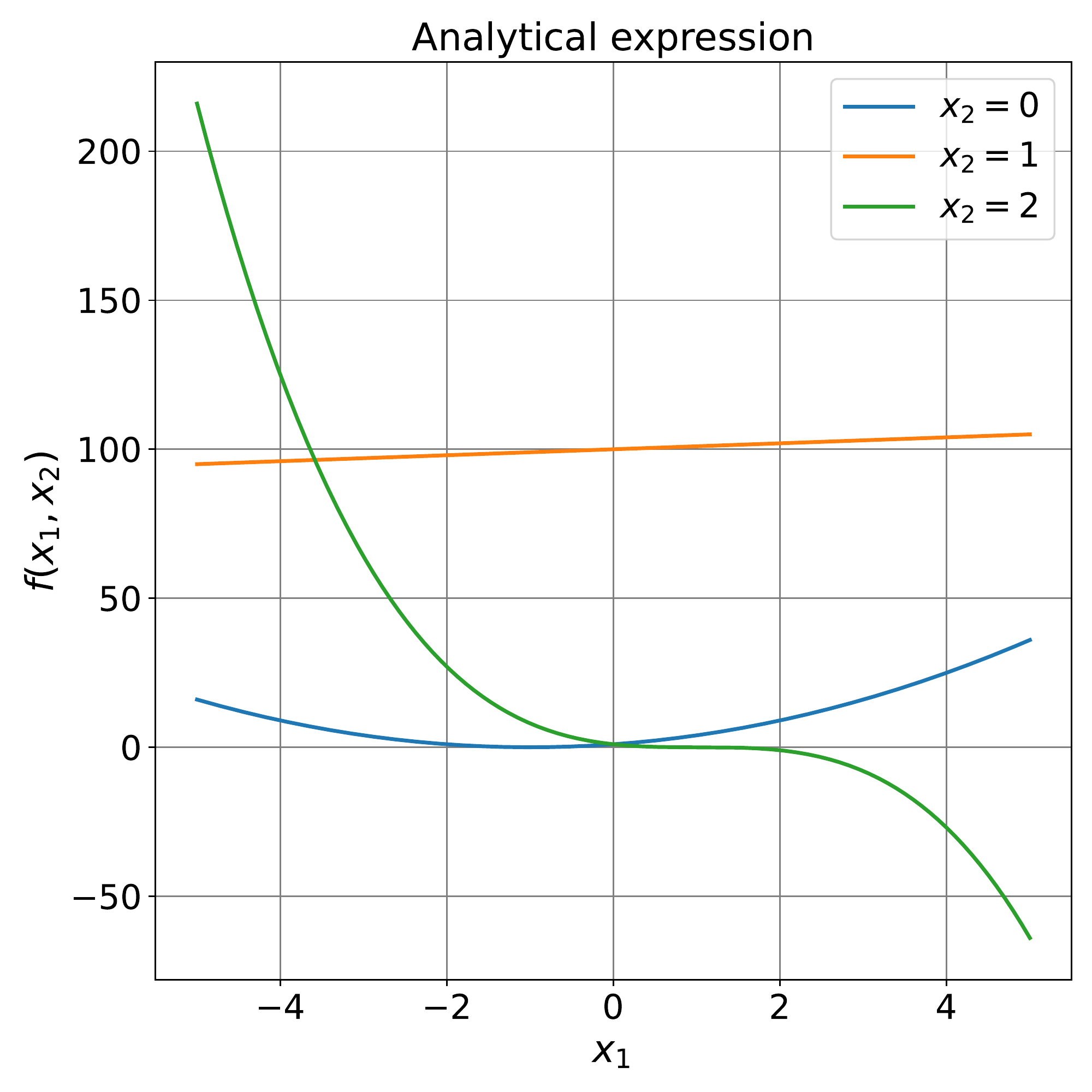}}\quad\quad 
    \subfigure[]{\includegraphics[width=0.45\textwidth]{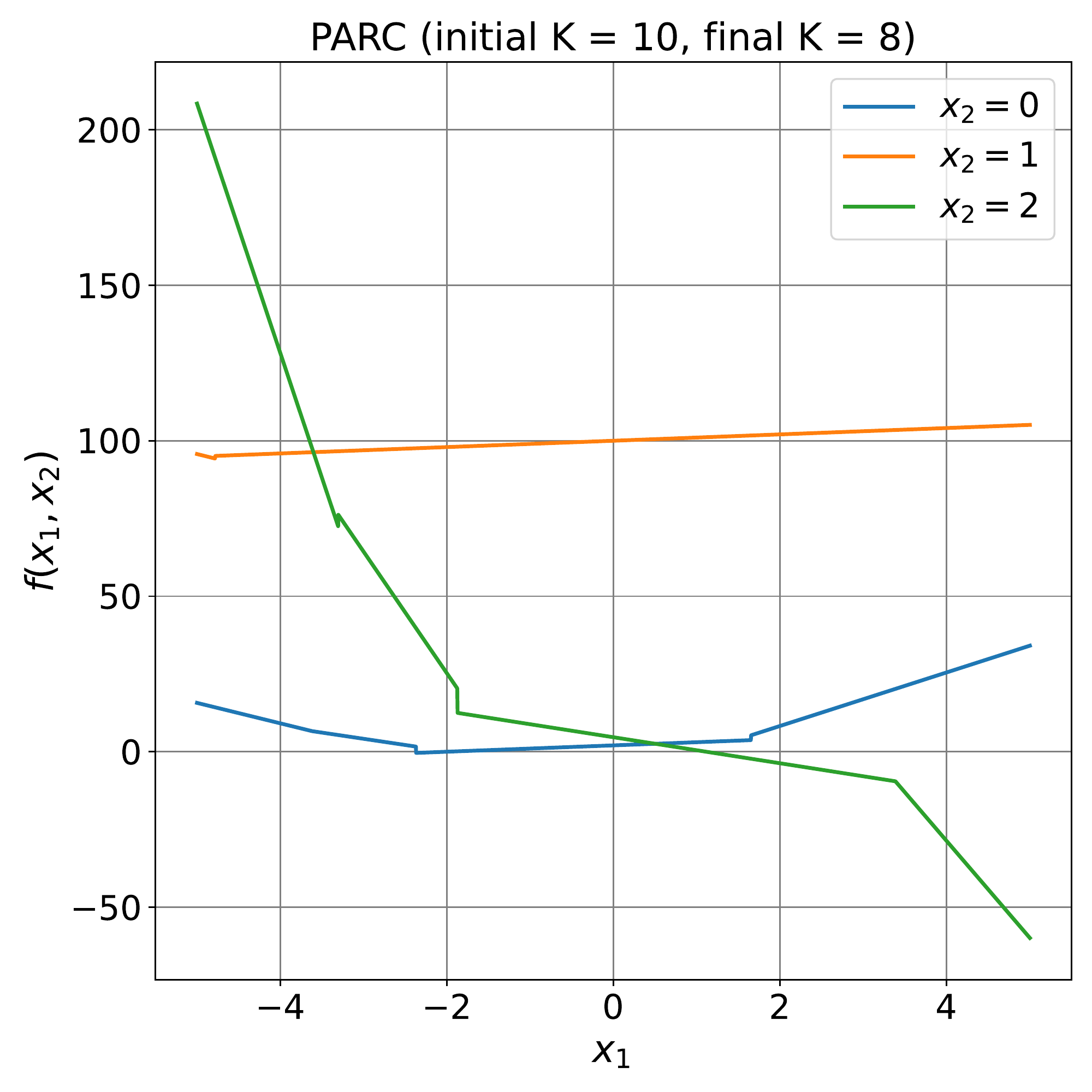}}
    \caption{(a) Function~\eqref{eq:syn_fun_illustrate} evaluated analytically, and (b) Function~\eqref{eq:syn_fun_illustrate} evaluated by the surrogates fitted by PARC with $K = 8$ final partitions. }
    \label{fig:pwe_fun_fitting_2}
\end{figure}

For the mixed continuous and categorical domain, we consider the following synthetic function:
\begin{equation}\label{eq:syn_fun_illustrate}
\begin{split}
    f(x_1,x_2) &=  \begin{cases}
 				x_1^2 + 2x_1 +1   & x_2= 0\\
				x_1 + 100 &x_2= 1\\
                    (1-x_1)^3   & x_2= 2
			\end{cases} \\
   -5 & \leq x_1 \leq 5, \ x_2\in \{0, 1, 2\},
\end{split}
\end{equation}
where with different categorical values of $x_2$, function evaluations ($f(x_1,x_2)$) can vary significantly. We use 960 randomly generated training samples with 10 initial partitions ($K$) to fit the PWA surrogates. Figure~\ref{fig:parc_partition_2} shows the final polyhedral partition induced by \eqref{eq:pwa_partition}. In this case, we observe that 2 initial partitions were discarded during surrogate fitting, resulting in 8 remaining partitions. In Figure~\ref{fig:pwe_fun_fitting_2}, we show the Function~\eqref{eq:syn_fun_illustrate} fitted analytically, and fitted by PARC with $K = 8$  partitions, where we observe that PARC can make good predictions at different values of the categorical variable ($x_2$), despite their distinct characteristics.

We also remark that 
the goal of the current study is to obtain a highly accurate approximation of the objective function \emph{around the global optimal solution} and not necessarily over the entire domain of $\bar X$, which usually requires much fewer samples. 
It is because, as the algorithm adaptively queries points to test from the domain, the partitions associated with higher function evaluations (less relevant for optimization) will be sampled less frequently, and accurate prediction models for these partitions are not necessary.  On the other hand, the regions with promising test points will be more frequently visited, resulting in better (more accurate) PWA surrogates within these partitions.  An illustrative example is shown in Appendix~\ref{sec:app_ill3} to demonstrate this remark.

\subsubsection{Mixed-integer Linear Encoding of the Surrogate}
\label{sec:MILP-surrogate}
After learning the coefficients of $\phi$ and $\hat f$ by applying the PARC algorithm, 
in order to optimize over the surrogate function to acquire a new sample $\bar X_{N+1}$
by MILP, as we will describe in Section~\ref{subsec:acq_fun}, we introduce $K$ binary variables $\zeta_j\in\{0,1\}$ and $K$ real variables $v_j\in\rr$, for $j=1,\ldots,K$. Here, $\zeta_j=1$ if and only if $\bar  X_{N+1}$ belongs to the $j$th polyhedral region of the partition induced by $\phi$. The PWA separation function $\phi$ can be modeled by the following mixed-integer inequalities via the big-M method~\cite{BemPARC21}:
\begin{equation}\label{eq:PWL_encode}
\begin{split}
&\omega_j\tr\bar X_{N+1} + \gamma_j\geq \omega_h\tr\bar X_{N+1} + \gamma_h - M_{\phi}(1-\zeta_j),\ \forall h=1,\ldots,K, h\neq j\\
&\sum_{j=1}^K\zeta_j=1,
\end{split}
\end{equation}
where $M_{\phi}$ is a large-enough constant, \emph{i.e.}, satisfies the inequality
\[
	M_{\phi}\geq \max_{j,h=1,\ldots,K,\ \bar X\in D} (\omega_h-\omega_j)\tr\bar X+ 
	\gamma_h - \gamma_j.
\] 
The PWA surrogate function $\hat f$ can be modeled by setting 
\[
	\hat{f}(\bar X) = \sum_{j=1}^K\zeta_j(a_j\tr\bar X_{N+1}+b_j) = \sum_{j=1}^K v_j,
\]
subject to
\begin{equation}\label{eq:surrogate_encode}%
\begin{split}
&v_j\leq a_j\tr\bar X_{N+1}+b_j-M_{sj}^-(1-\zeta_{j})\\
&v_j\geq a_j\tr\bar X_{N+1}+b_j-M_{sj}^+(1-\zeta_{j})\\
&v_j\geq M_{sj}^- \zeta_{j}\\
&v_j\leq M_{sj}^+ \zeta_{j},
\end{split}
\end{equation}
where $M_{sj}^+$, $M_{sj}^-$ are large-enough constants satisfying the inequalities
\[
	M_{sj}^+\geq \max_{\bar X\in D} a_j\tr\bar X+b_j,\quad 
	M_{sj}^-\leq \min_{\bar X\in D} a_j\tr\bar X+b_j,
\] 
for $j=1,\ldots,K$.

\subsection{Exploration Function}\label{subsec:exploration_fun}
Solely minimizing the surrogate function $\hat f$ may easily miss the global optimum.
In order to properly explore the admissible set $D$ we need to introduce
an \emph{exploration function} $E:\bar\Omega\mapsto\rr$. Due to the different numerical properties of continuous, integer, and categorical variables, we consider different exploration strategies for each of them that admit a MILP representation. Specifically, we use a distance-based exploration method for continuous and integer variables if the latter are not one-hot encoded (as described in Section~\ref{sec:scenario2}) and a frequency-based exploration method for one-hot encoded categorical and integer variables (in the alternative scenario described in Section~\ref{sec:scenario1}). In the following, we discuss the distance-based and frequency-based methods in a general manner, and we will dive into specifics of the exploration functions for our problem of interest when we discuss the acquisition function in Section~\ref{subsec:acq_fun}.

\subsubsection{Distance-based Exploration: ``Max-box'' Method} \label{subsec:max_box}
We want to define a function ${E_{ct}:\rr^{n_{ct}}\mapsto\rr}$ mapping a generic
numeric vector $\bar x\in\rr^{n_{ct}}$ into a nonnegative value $E_{ct}(\bar x)$ with the following features: ($i$) is zero at given samples $\bar x_1,\ldots,\bar x_N$; ($ii$) grows away from them; and ($iii$) admits a PWA representation. To this end, we consider the boxes $B_i(\beta_{ct})=\{\bar x:\ \|\bar x-\bar x_i\|_\infty\leq\beta_{ct}\}$
and set $E_{ct}(\bar x) = \min\{\beta_{ct}\geq 0:\ \bar x\in B_i(\beta_{ct})$
for some $i=1,\ldots,N\}$. Then, maximizing $E_{ct}(\bar x)$ is equivalent to finding the largest value $\beta_{ct}$ and a vector $\bar x^*$ outside the interior of all boxes $B_i(\beta_{ct})$, a problem that can be solved by the following MILP
\begin{equation}\label{eq:exploration_maxBox}%
\begin{array}{rcll}
	\bar x^* \in &\arg \  \underset{\bar x,\beta_{ct},\delta^+,\delta^-}{\max}  \ \beta_{ct}  \\[1em]
    \rm{s.t} &\quad \bar x^l-\bar x_i^l\geq \beta_{ct} - M_E(1-\delta_{il}^+),&\forall l=1,\dots, n_{ct},\quad \forall  i=1,\ldots,N\\[1em]
    &\quad -\bar x^l+\bar x_i^l\geq \beta_{ct} - M_E(1-\delta_{il}^-),&\forall l=1,\dots, n_{ct},\quad \forall  i=1,\ldots,N\\[1em]
    &\quad \delta_{il}^+\leq 1-\delta_{il}^-,&\forall l=1,\dots, n_{ct},\quad\forall  i=1,\ldots,N\\[0.5em]
        &\quad \displaystyle{\sum_{l=1}^{n_{ct}}\delta_{il}^++\delta_{il}^-\geq 1},  & \forall  i=1,\ldots,N\\[1.5em]
    &\quad \beta_{ct} \geq 0, \ \ \bar x \in D,
\end{array}
\end{equation}
where $l$ denotes the $l$th component of vector $\bar x$, 
$\delta_{il}^-,\delta_{il}^+\in\{0,1\}$
are auxiliary optimization variables introduced to model the violation
of at least one of the linear inequalities that define the box $B_i(\beta_{ct})$, 
and $M_E$ is a large-enough constant satisfying the following inequality 
\[
    M_E  \geq  2 \left( \underset{l = 1, \dots, n_{ct}}{\max} \bar u^l_x - \underset{l = 1, \dots, n_{ct}}{\min} \bar \ell^l_x \right),
\]
where $\bar u^l_x$ and $\bar \ell^l_x$ are the upper and lower bounds, respectively, of the $l$th component of vector $\bar x$.

\begin{figure}[hbt!]
    \centering
    \includegraphics[width=0.8\textwidth]{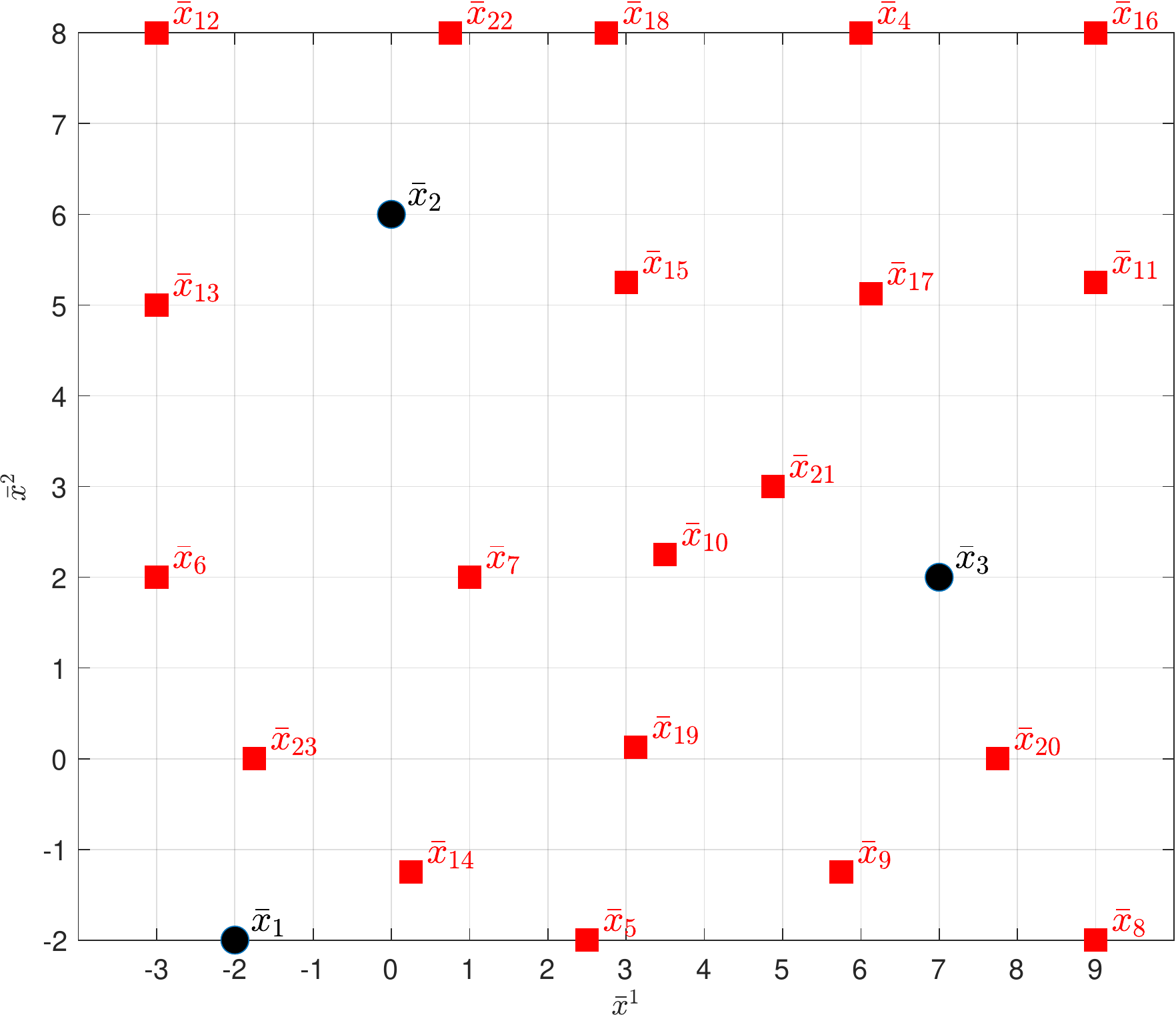}
    \caption{Illustrative example of the max-box exploration function in 2D. The black dots denote the initial samples. The red squares denote the samples generated using the max-box exploration method. The subscript number indicates the order of the point generated.}
    \label{fig:max_box_exp}
\end{figure}

Figure~\ref{fig:max_box_exp} shows an example, where we apply the max-box exploration method to $\bar x\in\rr^2$ and $D=[-3,9]\times[-2,8]$. We start with three existing samples $\bar x_1,\bar x_2,\ {\rm and}\ \bar x_3$. After 20 samples, we get the samples reported in the figure, which shows that, indeed, the max-box exploration method effectively explores the feasible region $D$.

\subsubsection{Frequency-based Exploration: ``Hamming distance'' Method}
Unlike the case of continuous variables treated in the previous section, to account for the frequency of occurrence of a particular combination of binary variables
we use the Hamming distance, defined as follows: given two binary vectors $z = [z^1,\dots,z^{d}]\tr\in\{0,1\}^d$ and $z_i = [z_i^1\ \ldots\ z_i^{d}]\tr\in\{0,1\}^d$, the Hamming distance between $z$ and $z_i$ is defined by
the number of different components between them
\begin{equation}\label{eq:hamm_dist}
d_H(z,z_i) = \sum_{m=1}^{d}\lvert z^m - z_i^m\rvert,
\end{equation}
which can be encoded as the following linear expression
\begin{equation} \label{eq:hamm_dist_encode}
d_H(z,z_i) = \sum_{m: z_i^m = 0} z^m + \sum_{m: z_i^m = 1} (1-z^m).
\end{equation}
We consider the exploration function $E_{dt}: \{0, 1\}^{d} \mapsto \rr $ 
such that $E_{dt}(z)$ quantifies the average number of different binary components 
between $z$ and the given $N$ vectors $z_1$, $\ldots$, $z_N$
\[
	E_{dt}(z)=\frac{1}{dN}\sum_{i=1}^Nd_H(z,z_i).
\]
Hence, a binary vector $z^*$ with maximum average Hamming distance $E_{dt}(z^*)$ from the
current samples $z_1, \dots, z_N$ can be determined by solving the following MILP
\begin{equation}\label{eq:exploration_HD}%
	z^* \in \arg \max_{z\in D}\  E_{dt}(z).
\end{equation}
Table~\ref{tab:illust_hamming} shows an example in which we have three categorical variables $Z = [Z^1, Z^2, Z^3]$, where $Z^1\in\{A,B\}$, $Z^2\in\{A,B,C,D,E\}$, and $Z^3\in\{A,B,C\}$. We start with three initial samples $Z_1 = [A, E, C]$, $Z_2 = [B, B, B]$, and $Z_3 =  [A,  D,  C]$. First, we binary encode the categorical variables, getting the corresponding vectors $z_1$, $z_2$, $z_3\in\{0,1\}^{10}$. Then, we solve the optimization problem~\eqref{eq:exploration_HD} to identify $z_4=z^*$ and its corresponding decoded form $Z_4$. The table shows the categorical values $Z_4,\ldots,Z_{23}$ generated in 20 subsequent sampling steps, which shows that a diverse set of categorical variables are obtained when applying the Hamming distance exploration method.

\begin{table}[hbt!]
\caption{Illustrative example of the Hamming distance exploration function.}
\centering
\begin{tabularx}{\linewidth}{>{\hsize=2\hsize}X *{10}{>{\hsize=0.7\hsize}X}}
\hline
\textbf{Sample} & \textbf{1}  & \textbf{2}  & \textbf{3}  & \textbf{4}  & \textbf{5}  & \textbf{6}  & \textbf{7}  & \textbf{8}  & \textbf{9}  & \textbf{10} \\ \hline
\textbf{$Z^1$}     & B           & A           & B           & A           & B           & A           & B           & A           & B           & A           \\
\textbf{$Z^2$}     & A           & C           & D           & A           & E           & B           & C           & D           & A           & E           \\
\textbf{$Z^3$}     & A           & B           & A           & C           & A           & B           & C           & A           & B           & C           \\ \hline
\textbf{Sample} & \textbf{11} & \textbf{12} & \textbf{13} & \textbf{14} & \textbf{15} & \textbf{16} & \textbf{17} & \textbf{18} & \textbf{19} & \textbf{20} \\ \hline
\textbf{$Z^1$}     & B           & A           & B           & A           & B           & A           & B           & A           & B           & A           \\
\textbf{$Z^2$}     & B           & C           & D           & A           & E           & B           & C           & D           & A           & E           \\
\textbf{$Z^3$}            & A           & B           & C           & A           & B           & C           & A           & B           & C           & A           \\ \hline
\end{tabularx}
\label{tab:illust_hamming}
\end{table}

We remark that both exploration functions (distance-based and frequency-based) are independent from the surrogate and does not require explicit uncertainty measurements, making it flexible to be integrated with different types of surrogates (\emph{e.g.}, polynomials).


\subsection{Acquisition Function}\label{subsec:acq_fun}
%
The surrogate and exploration functions defined in Sections~\ref{subsec:PWA_surrogate} and~\ref{subsec:exploration_fun} can be combined into the following 
\emph{acquisition problem}
\begin{subequations}
\begin{equation}
    \text{find}\ \bar X^*\in \arg \min_{\bar X \in D}\hat f(\bar X)-\delta (E_{ct}(\bar x) + E_{dt}([\bar y\tr\  z\tr]\tr)),
\label{eq:acq_detailed-a}
\end{equation}
when integer variables are treated as categorical as described in Section~\ref{sec:scenario1}, or into
\begin{equation}
    \text{find}\ \begin{bmatrix}\bar X^*\\y^*\end{bmatrix} \in \arg \min_{\bar X \in D, y\in[\ell_y,u_y]\cap \zz} \hat f(\bar X)-\delta (E_{ct}([\bar x\tr \ \bar y\tr]\tr) + E_{dt}(z)),
\label{eq:acq_detailed-b}
\end{equation}%
\label{eq:acq_detailed}%
\end{subequations}
when the integer vector $y$ is scaled as described in Section~\ref{sec:scenario2}. In~\eqref{eq:acq_detailed}, the nonnegative scalar $\delta$ is called the \emph{exploration parameter} and decides the tradeoff between exploiting the surrogate $\hat f(\bar X)$ and promoting the exploration of the feasible domain $D$. In the sequel, we will refer to the cost function $a:\bar\Omega\mapsto\rr$
in~\eqref{eq:acq_detailed-a} or~\eqref{eq:acq_detailed-b} as the \emph{acquisition function}. By construction, Problem~\eqref{eq:acq_detailed} can be solved by MILP. 
An optimal vector $\bar X^*$ of its solution, once scaled and decoded back, defines the next sample $X_{N+1}$ to query for the corresponding function value $f_{N+1}=f(X_{N+1})$. Note that $X_{N+1}$ satisfies all the constraints in~\eqref{eq:opt_prob_1}
since $\bar X^*\in D$.

The direct formulation~\eqref{eq:acq_detailed} can be further improved 
to ease the selection of $\delta$ and make the exploration more homogenous with respect to all types of variables (continuous, integer, or categorical). In fact, 
the formulation in~\eqref{eq:acq_detailed} has the following possible drawbacks:
\begin{enumerate}[label=(\roman*)]
\item The relative magnitude between $\hat f(\bar X)$ and $E(\bar X)$ is hard or impossible to estimate a priori, making the value for the exploration parameter $\delta$ hard to select.
\item By using the same exploration parameter $\delta$ for the exploration function of each type of variable, we implicitly assumed that the relative magnitude of each exploration function is comparable, which may not be the case.
\item When the integer variables $y$ are not one-hot encoded as described in Section~\ref{sec:scenario2}, the max-box exploration function is applied to the combined vector $[\bar x\tr\ \bar y\tr]\tr$ (see~\eqref{eq:acq_detailed-b}) and two problems can occur. Firstly, as shown in~\eqref{eq:opt_prob_2}, even though $\bar y_i$ is a continuous variable, 
because of the presence of the corresponding auxiliary integer variable $y_i$,
it can only be changed in discrete steps, unlike the remaining variables $\bar x_j$.
As a result, when finding the max box, one unit change in an integer variable can be reflected as a more significant change of its corresponding scaled variable $\bar y$, therefore promoting the exploration of directions with more variations in the integer variables $\bar y_i$ than in the continuous variables $\bar x_j$.
Secondly, due to possibly different lower and upper bounds and therefore scaling factors of integer variables, unit changes of them may cause different changes in size of their corresponding scaled variables.
\end{enumerate}
To address the aforementioned issues, given $N$ samples $(\bar X_i, f(S(\bar X_i))$, for
$i=1$, $\ldots$, $N$, we reformulate the acquisition problems~\eqref{eq:acq_detailed}, respectively, as follows:
\begin{subequations}
\begin{equation}
    \text{find}\ \bar X_{N+1}\in \arg \min_{\bar X \in D}
\frac{\hat f(\bar X)}{\Delta F} -
\delta_1 E_{ct}(\bar x) - \delta_2 E_{ct}(\bar y)- \delta_3 E_{dt}(z),
\label{eq:acq_solve-a}
\end{equation}
\begin{equation}
    \text{find}\ \begin{bmatrix}\bar X_{N+1}\\y_{N+1}\end{bmatrix} \in \arg \min_{\bar X \in D, y\in[\ell_y,u_y]\cap \zz} 
\frac{\hat f(\bar X)}{\Delta F} -
\delta_1 E_{ct}(\bar x) - \delta_2 E_{dt}(\bar y) - \delta_3 E_{dt}(z),
\label{eq:acq_solve-b}
\end{equation}%
\label{eq:acq_solve}%
\end{subequations}
where
\[
	\Delta F = \max\left\{\max_{i = 1, \dots, N}f(X_i) - \min_{i = 1, \dots, N}f(X_i), \epsilon_{\Delta F} \right\},
\]
and $\epsilon_{\Delta F} >0$ is a threshold to prevent division by zero.
The scaling factor $\Delta F$ eases the selection of the exploration parameters $\delta_1$, $\delta_2$, and $\delta_3$ by making the surrogate term comparable to the exploration terms (cf.~\cite{Bem20}). 

An alternative to solve the optimization problem~\eqref{eq:acq_solve} in one step
is to consider only one exploration term at a time, therefore solving
the problem in three consecutive steps (that will be referred to as the ``multi-step'' approach), where, at each step, 
the problem is only solved with respect to one variable type. The remaining variables 
are treated as constants at either the value associated with the current best vector $X_{N^*_{\rm curr}}$ or at the new value optimized during the multi-step operation. 
The advantage of serializing the optimization is that the
relative value of $\delta_1$, $\delta_2$, and $\delta_3$ is no longer relevant, 
and therefore we set $\delta_1=\delta_2=\delta_3=\delta$, where $\delta$
is the only tradeoff hyperparameter to choose.
Empirical findings indicate that a multi-step optimization approach for the acquisition function often yields superior results compared to a one-step approach. The one-step approach poses challenges in effectively tuning the exploration parameters $\delta_1$, $\delta_2$, and $\delta_3$, which can hinder performance. Additionally, these findings also motivate our heuristic for different treatments of integer variables as discussed in Section~\ref{subsec:pre-processing}. When optimize the acquisition function in multiple steps, if the number of possible integer combinations is relatively small compared to the maximum allowed number of black-box function evaluations, the integer variable combinations may be exhaustively enumerable within the allowed evaluations. Treating them categorically helps prevent premature convergence and ensures a more thorough exploration of the solution space.

A further heuristic is applied to restrict the number of binary variables used to encode the max-box exploration function~\eqref{eq:exploration_maxBox} and therefore limit them as the number $N$ of samples increases. Specifically, given
an upper bound $N_{\rm Emax}$ defined by the user depending on the computational power available, we only consider the most recent $N_S$ samples in the exploration function (we will use $N_S=20$ in our experiments) when $N n_c \geq N_{\rm Emax}$ or $N n_{\rm int} \geq N_{\rm Emax}$ (when integer variables are not one-hot encoded).
We note that the surrogate function is approximated using all the existing samples.
The rationale behind the heuristic is that, as the number $N$ of queried samples grows, the surrogate itself should already discourage the exploration around the older
samples not included in the exploration term, where the surrogate function, most likely, takes large values.

\subsection{Initial Sampling Strategies}\label{subsec:initial_sampl}
The values of the initial samples $X_1,\ldots,X_{N_{\rm init}}$ can significantly impact the final solution $X^*$ obtained after $N_{\rm max}$ steps. Moreover, one of the main motivations of the proposed method is its ability to handle mixed-integer constraints on the optimization variables. We propose different initial sampling strategies to efficiently acquire $N_{\rm init}$ scattered feasible samples depending on the constraints and types of optimization variables present in the problem:
\begin{enumerate}[label=(\roman*)]
	\item When only box constraints are present, we use the \emph{Latin Hypercube Sampling} (LHS)~\cite{mckay_lhs} method as in~\cite{Bem20}.
	\item 
	When both box constraints and linear equality and/or inequality constraints are present, we consider the following alternatives:
	\begin{itemize}
		\item If only continuous variables are present, we use the Double Description Method~\cite{motzkin1953double} to generate the $n_V$ vertices of the convex polytope given by the linear and box constraints. If $n_V<N_{\rm init}$ and only inequality constraints are involved in~\eqref{eq:opt_prob_1}, additional feasible samples can be generated via linear combinations of the vertices; if
		$n_V<N_{\rm init}$ and equality constraints are also present, the generated $n_V$ vertices can be used to define initial boxes and additional scattered feasible samples are generated by solving MILPs sequentially with the max-box exploration function discussed in Section~\ref{subsec:max_box} as the objective function. In this way, the constraints in~\eqref{eq:opt_prob_1} can be enforced in the formulation.
		\item If integer and/or categorical variables are also present, the algorithm first attempts to generate samples using LHS and filters out the infeasible ones. If the number of generated feasible samples is insufficient, which may happen when the constraints are hard to fulfill by random sampling, we generate scattered samples by solving MILPs sequentially using the exploration functions discussed in Section~\ref{subsec:exploration_fun} as the objective functions and incorporating the given constraints to ensure sample feasibility.
	\end{itemize}
\end{enumerate}

%

\section{Preference-based Learning}\label{sec:pref_learning}
We want to extend the global optimization method introduced in the previous sections to handle cases in which quantifying an objective function $f(X)$ as in~\eqref{eq:opt_prob_1} can be hard, if not impossible. For example, if multiple objectives are involved,
defining their relative weights a priori to form a single objective function can be difficult. In such cases,
expressing a \emph{preference} between the outcomes of two decision vectors $X_1$ and $X_2$ can be much simpler for a decision maker than quantifying the outcomes. Accordingly, we define the following \emph{preference function} $\pi:\Omega\times \Omega\to\{-1,0,1\}$ (cf.~\cite{BemPig20})
\begin{equation}
\pi(X_1,X_2)=\left\{\ba{ll}
-1 &\mbox{if $X_1$ ``better'' than $X_2$}\\
0 &\mbox{if $X_1$ ``as good as'' $X_2$}\\
1 &\mbox{if $X_2$ ``better'' than $X_1$}.
\ea\right.
\label{eq:pref_fun}
\end{equation}
In this case, we are interested in finding a feasible optimization vector $X^*$ that wins or ties the pairwise comparisons with any other feasible $X$ according to the preference function $\pi$, \emph{i.e.}, the optimization problem~\eqref{eq:opt_prob_1}
is replaced by
\begin{equation}
\text{find } X^*\ \mbox{such that}\ \pi(X^*,X)\leq 0,\ \forall X \in D.
\label{eq:glob-opt-pref}
\end{equation} 
We describe next a variant of Algorithm~\ref{algo:pwas}, that we call as PWASp, for solving the preference-based optimization problem~\eqref{eq:glob-opt-pref}.

Let the optimization vector $X$ be first pre-processed to $\bar X$ (\emph{e.g.}, scaling and/or encoding) as described in Section~\ref{subsec:pre-processing}. Given $N$ samples $\bar X_1, \dots, \bar X_N$ and $M_c$ preferences $\pi(S(\bar X_{1,k}),S(\bar X_{2,k})) \in\{-1,0,1\}$, for $k = 1,\dots, M_c$, where ${M_c = N -1}$, we aim to fit a PWA surrogate model reflecting the preference relations among different samples. Since function evaluations are not available, here, the preferences $\pi(S(\bar X_{1,k}),S(\bar X_{2,k}))$  are used to shape the surrogate function $\hat{f}(\bar X)$ by imposing the following constraints:
\begin{equation}\label{eq:pref_relation}
\begin{split}
\hat{f}(\bar X_{1,k}) \leq \hat{f}(\bar X_{2,k}) - \sigma &\quad \forall k:\ \pi(S(\bar X_{1,k}),S(\bar X_{2,k}))=-1\\
\hat{f}(\bar X_{2,k}) \leq \hat{f}(\bar X_{1,k}) - \sigma  &\quad \forall k:\ \pi(S(\bar X_{1,k}),S(\bar X_{2,k}))=1\\
|\hat{f}(\bar X_{1,k}) - \hat{f}(\bar X_{2,k})| \leq \sigma  &\quad \forall k:\ \pi(S(\bar X_{1,k}),S(\bar X_{2,k}))=0,		
\end{split}
\end{equation} 
where $(S(\bar X_{1,k}),S(\bar X_{2,k})) = (X_{1,k},X_{2,k})$ are pairs of compared samples, with $X_{1,k}$, $X_{2,k} \in\{X_1,\ldots,X_N\}$,
for $k=1,\ldots,M_c$. Here, $\sigma>0$ is a given constant, used to avoid the trivial solution $\hat f(\bar X)\equiv 0$. 

To identify the PWA separation function $\phi(\bar X)$, we first use K-means~\cite{lloyd1982least} to cluster the samples, and then use softmax regression~\cite{cox1951some,theil1969multinomial} to fit the coefficients. The assignment $j(\bar X)$ of each sample $\bar X$ to each region of the partition is then determined. Following that, different from the PARC algorithm, we determine the coefficients $a_j$, $b_j$ defining the PWA surrogate function $\hat{f}(\bar X)$ by minimizing the sum $\sum_{k=1}^{M_c}\epsilon_k$ of the violations of the preference constraints~\eqref{eq:pref_relation}
under an additional $\ell_\infty$-regularization term. Specifically, the coefficients
$a_j,b_j$ are obtained by solving the following linear programming (LP) problem:
\begin{equation}\label{eq:PWA_fit}
\begin{split}
   \min_{\epsilon_k,\xi, a, b} &\quad \sum_{k=1}^{M_c}\epsilon_k + \alpha \xi \\
   \text{s.t.}  &\quad \hat{f}(\bar X_{1,k}) + \sigma \leq \hat{f}(\bar X_{2,k}) + \epsilon_k \quad\  \forall k:\ \pi(X_{1,k},X_{2,k})=-1\\
  				 &\quad \hat{f}(\bar X_{2,k}) + \sigma \leq \hat{f}(\bar X_{1,k}) + \epsilon_k \quad \ \forall k:\ \pi(X_{1,k},X_{2,k})=1\\
  				  &\quad |\hat{f}(\bar X_{1,k}) - \hat{f}(\bar X_{2,k})| \leq \sigma + \epsilon_k \quad \forall k:\ \pi(X_{1,k},X_{2,k})=0\\				 
   				&\quad  \xi \geq\pm a_{j}^l,\ l=1,\ldots,n\\
    			&\quad	\xi \geq\pm b_{j},
\end{split}
\end{equation} 
where $\alpha >0$ is the regularization parameter; $\xi \in \rr$ is a new optimization variable introduced to linearly encode the $\ell_\infty$-regularization of the coefficients; and $l$ denotes the $l$th component of the vector.

Once the surrogate model is obtained, the same procedure as in PWAS can be followed to construct the acquisition function. This acquisition function is then optimized to identify the next sample, $X_{N+1} = S(\bar X_{N+1})$. The new sample is subsequently compared with the current best vector, $X_{N^*_{\rm curr}}$, to obtain the new preference assessment. The various steps involved in PWASp are summarized in Algorithm~\ref{algo:pwasp}.

\begin{algorithm}[hb!]
	\caption{PWASp: Preference-based Optimization Using Piecewise Affine Surrogates}
	\label{algo:pwasp}
	~~\textbf{Input}: Lower and upper bounds $\ell_x,u_x,\ell_y,u_y$; linear constraint matrices $A_{\rm eq}$, $B_{\rm eq}$, $C_{\rm eq}$, $A_{\rm ineq}$, $B_{\rm ineq}$, and $C_{\rm ineq}$; linear constraint right-hand-side vectors $b_{\rm eq}$ and $b_{\rm ineq}$; number $n_d$ of categorical variables and $n_i$ of possible categories, for $i=1,\ldots,n_d$; initial number $K$ of polyhedral partitions;
	number $N_{\rm init}\geq 2$ of initial samples to compare; maximum number $N_{\rm max}-1$ of comparisons, $N_{\rm max}\geq N_{\rm init}$; $\delta_1 \geq 0$, $\delta_2 \geq 0$ and $\delta_3 \geq 0$ if solve~\eqref{eq:acq_solve} in one step or $\delta \geq 0$ if solve~\eqref{eq:acq_solve} in multiple steps; solving strategy for~\eqref{eq:acq_solve}: \{``one-step'' or ``multi-steps''\}.
	\vspace*{.1cm}\hrule\vspace*{.1cm}
	\begin{enumerate}
		\item Pre-process the optimization variables as described in Section~\ref{subsec:pre-processing};
		\item  $N\leftarrow 1$, $i^*\leftarrow 1$;
		\item \label{algo:initial} Generate $N_{\rm init}$ random and encoded samples $\bar X=\{\bar X_1,\ldots, \bar X_{N_{\rm init}}\}$ using one of the initial sampling methods described in Section~\ref{subsec:initial_sampl} based on the problem setup;
			\item \textbf{While} $N < N_{\rm max}$ \textbf{do}
			\begin{enumerate}
				\item \textbf{If} $N \geq N_{\rm init}$ \textbf{then}
				\begin{enumerate}
					\item Update and fit the PWA separation function $\phi$ and PWA surrogate function $\hat f$ as described in Section~\ref{sec:pref_learning};
					\item Define the acquisition function $a$ as in~\eqref{eq:acq_solve};
					\item Solve the MILP problem~\eqref{eq:acq_solve} and get $\bar X_{N+1}$, either in one step or multiple steps;
				\end{enumerate}
				\item $i(N) \leftarrow i^*$, $j(N) \leftarrow N+1$
				\item Query preference $\pi(X_{i(N)}, X_{j(N)})$; 
				\item \textbf{If} $\pi(X_{i(N)}, X_{j(N)}) = 1$ \textbf{then set} $i^* \leftarrow j(N)$ 
				\item $N\leftarrow N + 1$;
			\end{enumerate}
		\item \textbf{End}.
	\end{enumerate}
	\vspace*{.1cm}\hrule\vspace*{.1cm}
	~~\textbf{Output}: Best vector $X^*=X_{i^*}$ encountered.
\end{algorithm}

\section{Optimization Benchmarks}\label{sec:benchmark}
To illustrate the effectiveness of PWAS and PWASp in solving the target problems~\eqref{eq:opt_prob_1} and~\eqref{eq:glob-opt-pref}, we have considered different mixed-variable global optimization benchmarks, including 
three unconstrained synthetic benchmarks as well as two unconstrained real-world benchmarks (taken from~\cite{ru2020bayesian}) and two constrained mixed-variable synthetic problems. 
Computations are performed on an Intel i7-8550U 1.8-GHz CPU laptop with 24GB of RAM. The MILP problem in the acquisition step is formulated with the PuLP library~\cite{mitchell2011pulp} and solved by Gurobi's MILP solver~\cite{gurobi}.

For each benchmark, the function evaluations are fed into PWAS to fit the surrogate, while the explicit function expressions remain unknown to PWAS.  As in~\cite{ru2020bayesian} the benchmark problems are solved via maximization, we use 
the values $-f(X)$ when running PWAS.  As for PWASp, the objective function serves as a synthetic decision-maker whose evaluations are only used to express the preference between two decision vectors, namely $\pi(X_1,X_2)=-1$
if $f(X_1)>f(X_2)$, $\pi(X_1,X_2)=1$ if $f(X_1)<f(X_2)$, or zero otherwise.  In other words, PWASp only has access to the queried preferences~\eqref{eq:pref_fun} and not the explicit function expressions nor their evaluations $f(X_N)$. Specifically, $N_{\rm max} -1$ pairwise comparisons are obtained for each benchmark when solved by PWASp. 

\begin{table}[b!]
\caption{Benchmark problem specifications.}
\centering
\begin{tabularx}{\linewidth}{>{\hsize=1.8\hsize}X *{3}{>{\hsize=0.5\hsize}X}{>{\hsize=1.2\hsize}X}}
\hline
\textbf{Benchmark}   & \textbf{$n_c$} & \textbf{$n_{\rm int}$} & \textbf{$n_d$} & \textbf{$n_i$}         \\ \hline
\textsf{Func-2C}                        & 2              & 0                      & 2              & \{3, 3\}               \\
\textsf{Func-3C}                        & 2              & 0                      & 3              & \{3, 3, 3\}            \\
\textsf{Ackley-5C}                      & 1              & 0                      & 5              & \{17, 17, 17, 17, 17\} \\
\textsf{XG-MNIST}                       & 4              & 1                      & 3              & \{2, 2, 2\}            \\
\textsf{NAS-CIFAR10}                    & 21             & 1                      & 5              & \{3, 3, 3\}            \\
\textsf{Horst6-hs044-modified}          & 3              & 4                      & 2              & \{3, 2\}               \\
\textsf{ros-cam-modified}               & 2              & 1                      & 2              & \{2, 2\}               \\ \hline
\end{tabularx}
\label{tab:benchmark_spec}
\end{table}

The performance of PWAS and PWASp on the unconstrained benchmarks are compared with the following solvers: CoCABO-auto~\cite{ru2020bayesian}, CoCABO-0.5~\cite{ru2020bayesian}, One-hot BO~\cite{gpyopt2016}, SMAC~\cite{hutter2011sequential}, TPE~\cite{bergstra2011algorithms}, and EXP3BO~\cite{gopakumar2018algorithmic} as noted in~\cite{ru2020bayesian} as well as MISO~\cite{muller2016miso} and NOMAD~\cite{abramson2011nomad,audet2021nomad}. A brief summary of each algorithm is provide in Appendix~\ref{sec:app_algorithm}. CoCABO-auto and CoCABO-0.5 are selected because the authors noted that they consistently show competitive performance~\cite{ru2020bayesian}. MISO and NOMAD are selected since they are noted as the best performers among all the algorithms tested in~\cite{ploskas2022review}. The settings of the first six algorithms compared are available in~\cite{ru2020bayesian}. As for MISO and NOMAD, we kept their default algorithm settings, with integer and categorical variables declared as referenced by the corresponding algorithms~\cite{muller2016miso,abramson2011nomad,audet2021nomad}.
To have fair comparisons, we performed the same initial and maximum number of black-box function evaluations ($N_{\rm init} = 20$ and $N_{\rm max} = 100$) as indicated in~\cite{ru2020bayesian} on the synthetic and real-world benchmarks for all the tested algorithms except for NOMAD. NOMAD starts the optimization process with an initial guess~\cite{audet2021nomad}. 
When solving the problems using PWAS or PWASp, the multi-step solution strategy is applied in the acquisition step with $\delta_1 = \delta_2 = \delta_3 = 0.05$ or $\delta_1 = \delta_2 = \delta_3 = 1$, respectively for PWAS and PWASp. The initial number $K$ of polyhedral partitions is set to $20$ for both PWAS and PWASp in all the benchmarks.
We stress that here the function evaluations for PWASp are solely reported for performance comparisons and are not attainable to PWASp during optimization.
Table~\ref{tab:benchmark_spec} summarizes the tested benchmark problems, while a detailed description of the benchmarks is reported in Appendix~\ref{sec:app_bench}. 

The optimization results obtained by CoCaBO-auto, CoCaBO-0.5, One-hot BO, SMAC and TPE for unconstrained benchmark were read from Figure 4 in~\cite{ru2020bayesian} using GetData Graph Digitizer~\cite{getData}. Regarding MISO and NOMAD (version 4), we retrieved their packages from the GitHub repository. We performed 20 random repetitions for the unconstrained synthetic problems (\textsf{Func-2C}, \textsf{Func-3C}, and \textsf{Ackley-5C}) and 10 random repetitions for the unconstrained real-world problems (\textsf{XG-MNIST} and \textsf{NAS-CIFAR10}) as reported in~\cite{ru2020bayesian}.

Regarding problems with constraints, we consider the two mixed-variable synthetic problems reported in Appendix~\ref{subsec:conts_synthetic}, whose specifications are also reported in Table~\ref{tab:benchmark_spec}. We do not consider other solvers than PWAS and PWASp as they either do not support constraint handling with mixed variables or allow infeasible samples during the optimization process. Thus, a systematic comparison is not performed for the constrained problems. Instead, the results are compared against the analytic global optimum.
Here, we set $N_{\rm max} = 100$ and $N_{\rm init} = \lceil N_{\rm max}/4 \rceil = 25$, $K = 20$ initial clusters, and the exploration parameters $\delta_1 = \delta_2 = \delta_3 = 0.05$ when using PWAS or $\delta_1 = \delta_2 = \delta_3 = 1$ with PWASp. We run PWAS and PWASp 20 times from different random seeds on these two problems. 

\subsection{Results and Discussions}\label{subsec:res_dis}
\begin{table}[tb!]
\caption{Best value found on benchmark \textsf{Func-2C}~\cite{ru2020bayesian} after 100 and 200 black-box function evaluations (maximum = 0.2063).}
\centering
\begin{tabularx}{\linewidth}{>{\hsize=1.4\hsize}X *{4}{>{\hsize=0.9\hsize}X}}
\hline
\multirow{2}{*}{Algorithm}  & \multicolumn{2}{l}{{After 100 evaluations}}   & \multicolumn{2}{l}{{After 200 evaluations}}  \\ \cline{2-5} 
                                     & mean           & std           & mean           & std          \\ \hline
\textsf{PWAS}                        & 0.2049                   & 0.0022                  &    0.2061                       &     0.0002321                   \\
\textsf{PWASp}                       & 0.1813                   & 0.0443                  & 0.1889                   & 0.0449                 \\
\textsf{CoCaBO-auto}                 & 0.1219                    & 0.0172                  & 0.2041                    & 0.0057                 \\
\textsf{CoCaBO-0.5}                  & 0.1352                    & 0.01620                 & 0.2041                    & 0.0057                 \\
\textsf{One-hot BO}                  & 0.009524                  & 0.02158                 & 0.01524                   & 0.02064                \\
\textsf{SMAC}                        & 0.06381                   & 0.01746                 & 0.07714                   & 0.01556                \\
\textsf{TPE}                         & 0.1273                    & 0.0184                  & 0.1743                    & 0.01650                \\
\textsf{EXP3BO}                      & 0.05524                   & 0.01429                 & 0.1105                    & 0.01650                \\ 
\textsf{MISO}                      &     0.2063               &    0.0000              &    0.2063                 &   0.0000              \\ 
\textsf{NOMAD}                      &   0.1700                 &    0.0736              &     0.1754                &     0.07557            \\ 
\hline
\end{tabularx}
\label{tab:fun_2c}

\bigskip

\caption{Best value found on benchmark \textsf{Func-3C}~\cite{ru2020bayesian} after 100 and 200 black-box function evaluations (maximum = 0.7221).}
\centering
\begin{tabularx}{\linewidth}{>{\hsize=1.4\hsize}X *{4}{>{\hsize=0.9\hsize}X}}\hline
\multirow{2}{*}{Algorithm}  & \multicolumn{2}{l}{{After 100 evaluations}}   & \multicolumn{2}{l}{{After 200 evaluations}}  \\ \cline{2-5} 
                                     & mean           & std           & mean           & std          \\ \hline
\textsf{PWAS}                        & 0.5282                    & 0.2117                  & 0.6450                    & 0.0972                 \\
\textsf{PWASp}                       & 0.4542                    & 0.2078                  & 0.5106                    & 0.1665                 \\
\textsf{CoCaBO-auto}                 & 0.4993                    & 0.0299                  & 0.6912                    & 0.0169                 \\
\textsf{CoCaBO-0.5}                  & 0.5371                    & 0.0503                  & 0.6991                    & 0.0205                 \\
\textsf{One-hot BO}                  & 0.007670                  & 0.04956                 & 0.1076                    & 0.0606                 \\
\textsf{SMAC}                        & 0.1084                    & 0.04016                 & 0.1965                    & 0.0339                 \\
\textsf{TPE}                         & 0.2672                    & 0.0472                  & 0.4914                    & 0.5308                 \\
\textsf{EXP3BO}                      & 0.1784                    & 0.0393                  & 0.2515                    & 0.0330                  \\  
\textsf{MISO}                      &      0.7221               &    0.0000              &   0.7221                   &   0.0000              \\
\textsf{NOMAD}                      &   0.6618                 &      0.1610            &      0.6860               &    0.1615             \\ 
\hline
\end{tabularx}
\label{tab:fun_3c}

\bigskip

\caption{Best value found on benchmark \textsf{Ackley-5C}~\cite{ru2020bayesian} after 100 and 200 black-box function evaluations (maximum = 0).}
\centering
\begin{tabularx}{\linewidth}{>{\hsize=1.4\hsize}X *{4}{>{\hsize=0.9\hsize}X}}
\hline
\multirow{2}{*}{Algorithm}  & \multicolumn{2}{l}{{After 100 evaluations}}   & \multicolumn{2}{l}{{After 200 evaluations}}  \\ \cline{2-5} 
                                     & mean           & std           & mean           & std          \\ \hline
\textsf{PWAS}                        & -1.1148                   & 0.4077                  &  -0.7108	                         &    0.3320                    \\
\textsf{PWASp}                       & -1.8857                   & 0.5795                  &    -1.6462	                      &  0.5422                       \\
\textsf{CoCaBO-auto}                 & -2.5120                   & 0.602                   & -1.9244                   & 0.5512                 \\
\textsf{CoCaBO-0.5}                  & -2.8415                   & 0.0488                  & -2.0073                   & 0.0488                 \\
\textsf{One-hot BO}                  & -3.076                    & 0.0483                  & -2.5341                   & 0.3024                 \\
\textsf{SMAC}                        & -3.0073                   & 0.2488                  & -1.710                    & 0.2393                 \\
\textsf{TPE}                         & -3.4659                   & 0.2000                     & -2.7976                   & 0.2487                 \\  
\textsf{MISO}                      &    -1.6389                &     0.1388             &     -1.5582                &    0.06218             \\ 
\textsf{NOMAD}                      &    -2.0175                &     0.2015             &     -1.5467                &    0.01437             \\ 
\hline
\multicolumn{5}{l}{\fontsize{8}{9}\selectfont Note: the reported values for CoCaBO-auto, CoCaBO-0.5, One-hot BO, SMAC, TPE  }\\
\multicolumn{5}{l}{\fontsize{8}{9}\selectfont  and EXP3BO are read from Figure 4 in~\cite{ru2020bayesian} using GetData Graph Digitizer~\cite{getData}.   }\\
\multicolumn{5}{l}{\fontsize{8}{9}\selectfont  Statistics are obtained over 20 runs.}\\
\end{tabularx}
\label{tab:Ackley-5C}
\end{table}

\begin{figure}[tbh!]
    \centering
    \includegraphics[width=.95\textwidth]{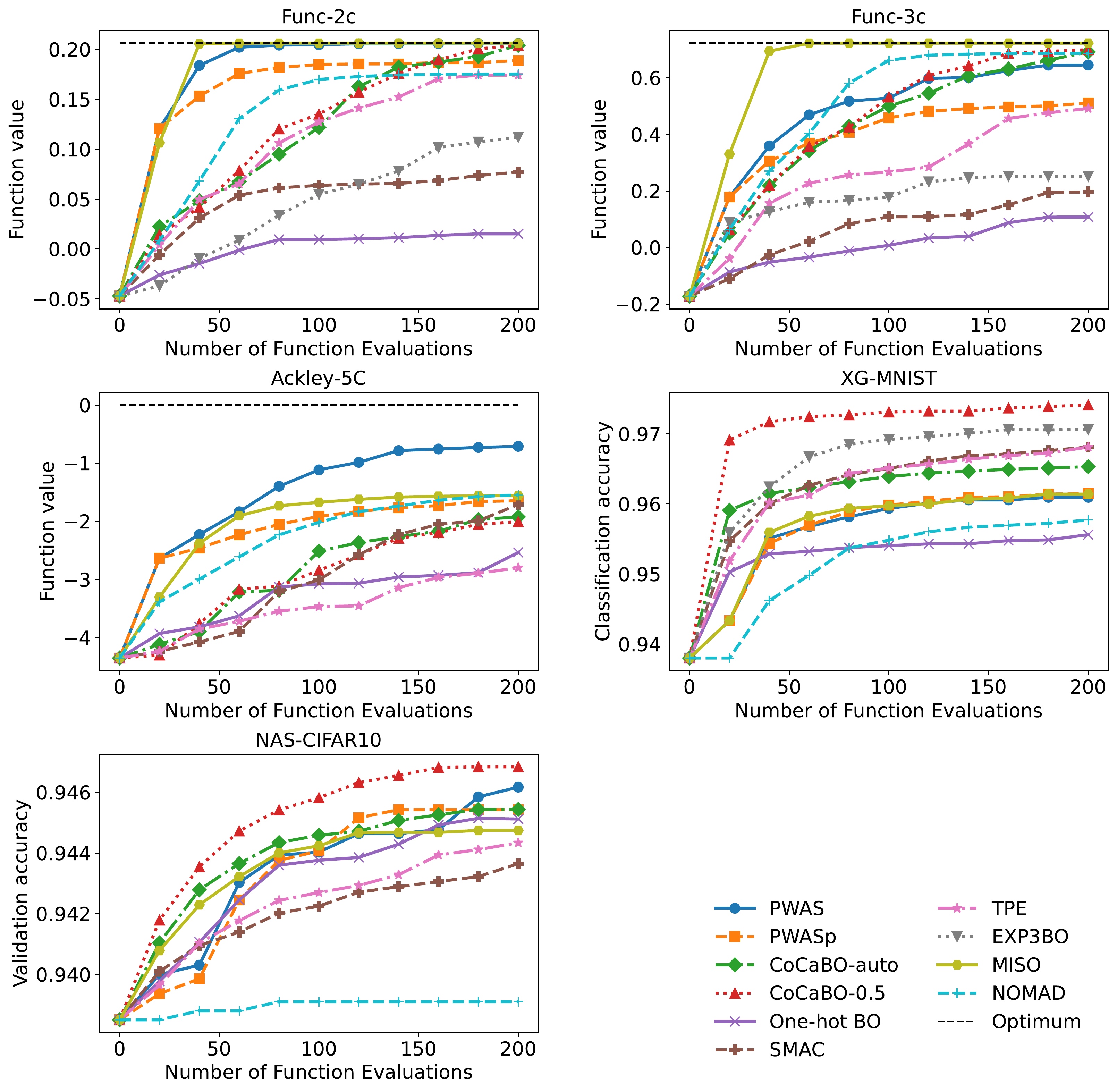}\quad\quad 

    \caption{Convergence graph of tested algorithms on the unconstrained benchmarks. The mean values are plotted, which are averaged over 20 runs for \textsf{Func-2C}, \textsf{Func-3C}, and \textsf{Ackley-5C}, and over 10 runs for \textsf{XG-MNIST} and \textsf{NAS-CIFAR10}.}
    \label{fig:convergence_plot}
\end{figure}

\begin{figure}[tbh!]
    \centering
    \includegraphics[width=.95\textwidth]{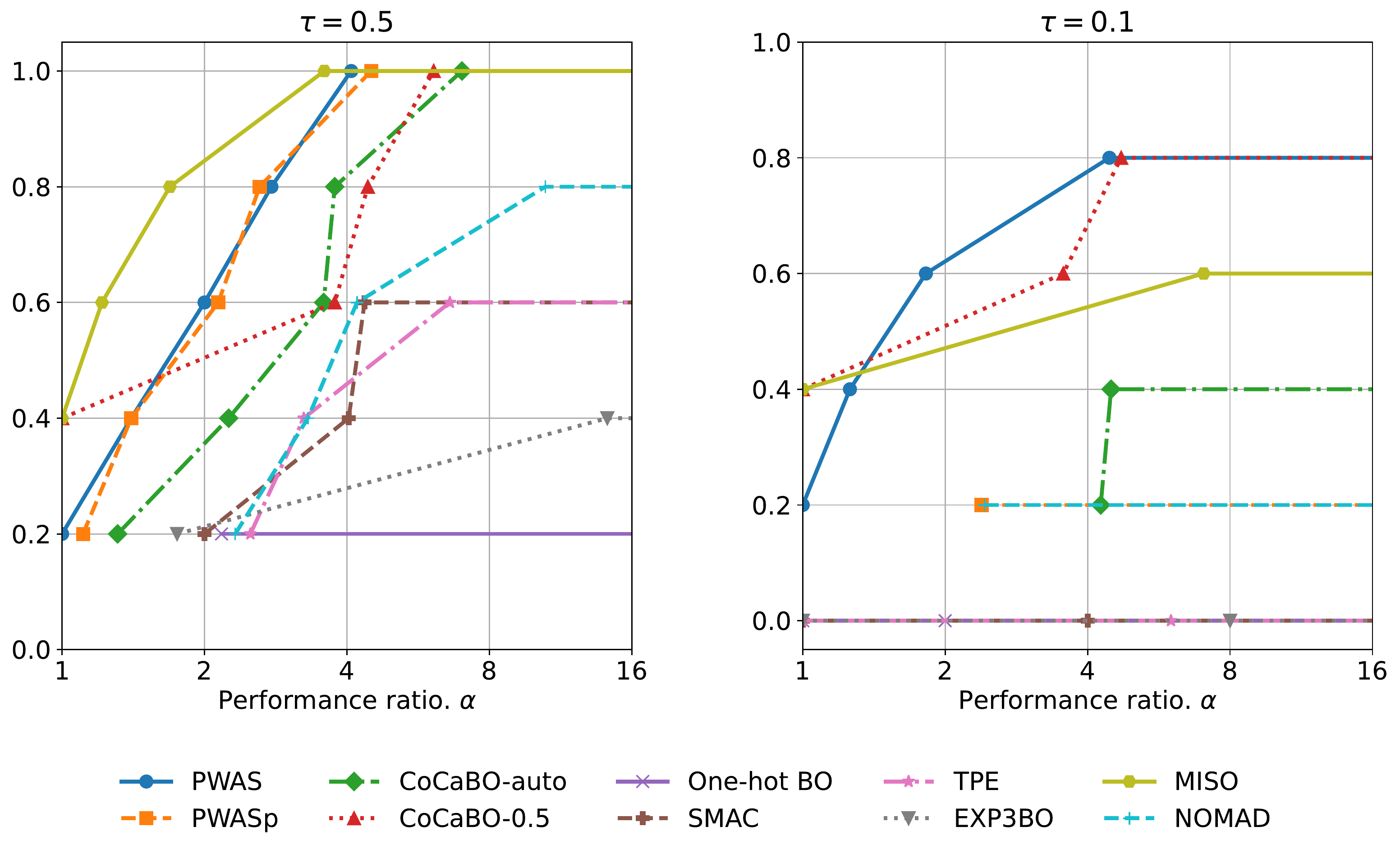}\quad\quad 

    \caption{Performance profiles of tested algorithms across varying performance ratios ($\alpha$) with convergence tolerance $\tau$ equals to 0.5 and 0.1. Performance ratio ($x$-axis) is in logarithmic scale with base 2. }
    \label{fig:performance_profile}
\end{figure}

\begin{figure}[tbh!]
    \centering
    \includegraphics[width=.95\textwidth]{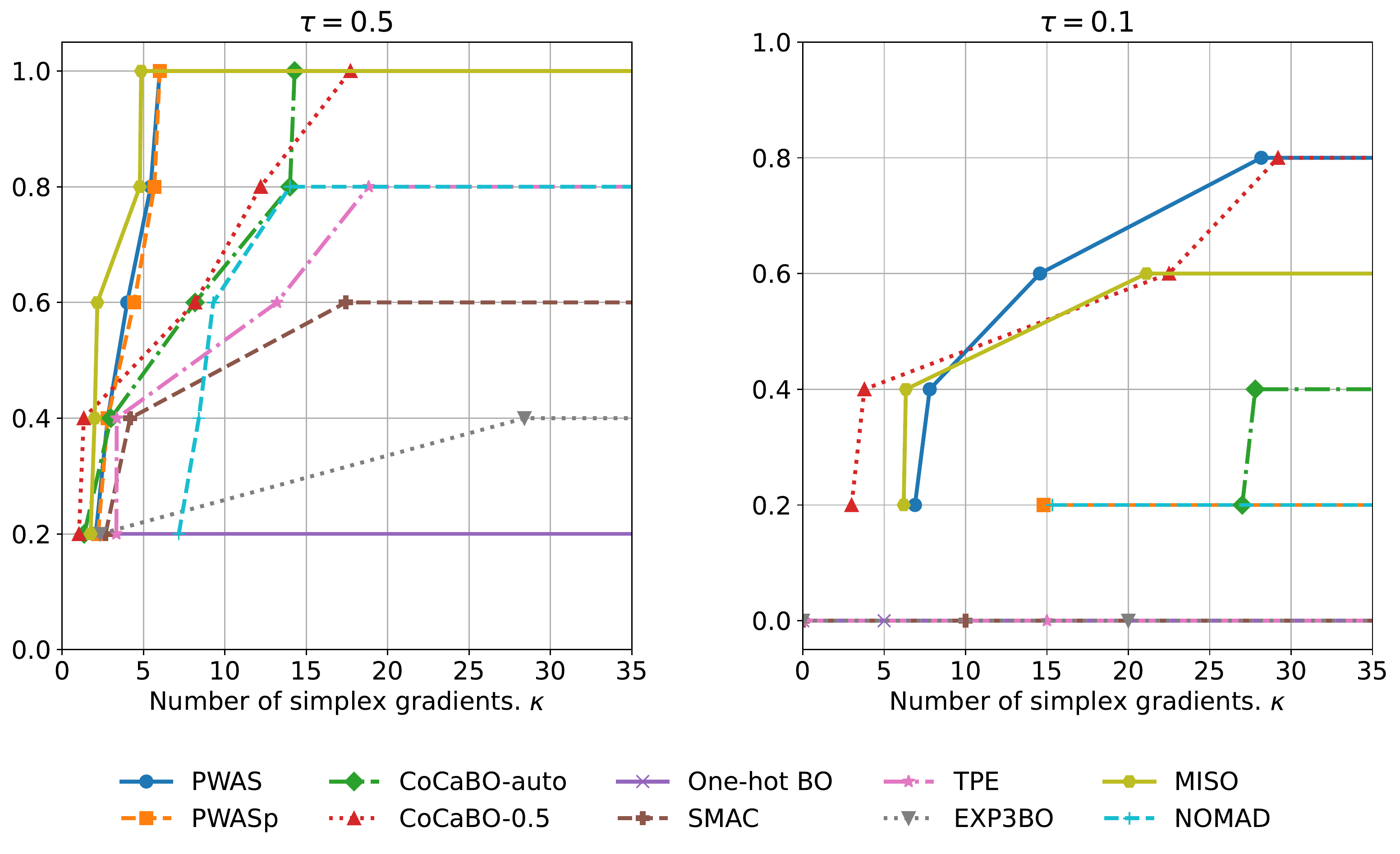}\quad\quad 

    \caption{Data profiles of tested algorithms as a function of the number of equivalent simplex gradients.}
    \label{fig:data_profile}
\end{figure}

The optimization results for unconstrained benchmarks are reported in Tables~\ref{tab:fun_2c}--\ref{tab:NAS-CIFAR10}, and illustrated in Figures~\ref{fig:convergence_plot}--\ref{fig:data_profile}.
The convergence plots (Figure~\ref{fig:convergence_plot}) show how different algorithms perform across various benchmarks in terms of increasing (reducing) the objective function value over time for maximization (minimization) problems. To avoid clutter, only mean values are plotted in Figure~\ref{fig:convergence_plot} with the standard deviations at evaluations 100 and 200 reported in Table~\ref{tab:fun_2c}--\ref{tab:NAS-CIFAR10}.
The performance and data profiles (Figure~\ref{fig:performance_profile} and~\ref{fig:data_profile}) are generated based on the guidelines noted in~\cite{more2009benchmarking} with convergence tolerance $\tau$ set to 0.5 and 0.1. 
Performance profiles compare the efficiency of algorithms by looking at the performance ratio (relative performance) across a set of problems, while 
data profiles illustrate the percentage of problems that can be solved as a function of the equivalent number of simplex gradients (function evaluations)~\cite{more2009benchmarking}. These profiles are particularly useful for comparing and analyzing the short-term behavior of algorithms, especially when computational resources are limited.

\begin{table}[b!]

\caption{Best value found on benchmark \textsf{XG-MNIST}~\cite{ru2020bayesian} after 100 and 200 black-box function evaluations.}
\centering
\begin{tabularx}{\linewidth}{>{\hsize=1.4\hsize}X *{4}{>{\hsize=0.9\hsize}X}} \hline
\multirow{2}{*}{Algorithm}  & \multicolumn{2}{l}{{After 100 evaluations}}   & \multicolumn{2}{l}{{After 200 evaluations}}  \\ \cline{2-5} 
                                     & mean           & std           & mean           & std          \\ \hline
\textsf{PWAS}                        & 0.9585                    & 0.0030                  &    0.9609                       &   0.0029                     \\
\textsf{PWASp}                       & 0.9576                    & 0.0036                  &    0.9615                      &      0.0028                  \\
\textsf{CoCaBO-auto}                 & 0.9639                    & 0.0004                  & 0.9653                    & 0.0004                 \\
\textsf{CoCaBO-0.5}                  & 0.9731                    & 0.0008                  & 0.9741                    & 0.0008                 \\
\textsf{One-hot BO}                  & 0.9541                    & 0.0019                  & 0.9556                    & 0.0015                 \\
\textsf{SMAC}                        & 0.9651                    & 0.0012                  & 0.9681                    & 0.0012                 \\
\textsf{TPE}                         & 0.9656                    & 0.0007                  & 0.9679                    & 0.0007                 \\
\textsf{EXP3BO}                      & 0.9691                    & 0.0005                  & 0.9706                    & 0.0005                 \\  
\textsf{MISO}                      &   0.9574                 &    0.0071              & 0.9594                     &         0.0078        \\ 
\textsf{NOMAD}                      &   0.9528                 &   0.0138               &         0.9564            &      0.0146           \\ 
\hline
\end{tabularx}
\label{tab:XG-MNIST}

\bigskip

\caption{Best value found on benchmark \textsf{NAS-CIFAR10}~\cite{ru2020bayesian} after 100 and 200 black-box function evaluations.}
\centering
\begin{tabularx}{\linewidth}{>{\hsize=1.4\hsize}X *{4}{>{\hsize=0.9\hsize}X}}
\hline
\multirow{2}{*}{Algorithm}  & \multicolumn{2}{l}{{After 100 evaluations}}   & \multicolumn{2}{l}{{After 200 evaluations}}  \\ \cline{2-5} 
                                     & mean           & std           & mean           & std          \\ \hline
\textsf{PWAS}        & 0.9440 & 0.0024 & 0.9462 & 0.0016 \\
\textsf{PWASp}       & 0.9409 & 0.0052 & 0.9454 & 0.0019 \\
\textsf{CoCaBO-auto} & 0.9446 & 0.0017 & 0.9454 & 0.0017 \\
\textsf{CoCaBO-0.5}  & 0.9458 & 0.0014 & 0.9468 & 0.0004 \\
\textsf{One-hot BO}  & 0.9438 & 0.0006 & 0.9451 & 0.0006 \\
\textsf{SMAC}        & 0.9422 & 0.0004 & 0.9436 & 0.0004 \\
\textsf{TPE}         & 0.9427 & 0.0006 & 0.9443 & 0.0007 \\  
\textsf{MISO}                                    &    0.9442                 &      0.0020     &   0.9447                 &    0.0035      \\ 
\textsf{NOMAD}                      &   0.9385                 &      0.0231            &     0.9391                &       0.0341         \\ 
\hline
\multicolumn{5}{l}{\fontsize{8}{9}\selectfont Note: the reported values for CoCaBO-auto, CoCaBO-0.5, One-hot BO, SMAC, TPE  }\\
\multicolumn{5}{l}{\fontsize{8}{9}\selectfont  and EXP3BO are read from Figure 4 in~\cite{ru2020bayesian} using GetData Graph Digitizer~\cite{getData}.  }\\
\multicolumn{5}{l}{\fontsize{8}{9}\selectfont  Statistics are obtained over 10 runs.}\\
\end{tabularx}
\label{tab:NAS-CIFAR10}

\end{table}


In general, PWAS and PWASp can effectively increase (decrease for minimization problems) the objective function values within a small number of function evaluations
or comparisons. In fact, as shown in Tables~\ref{tab:fun_2c}--\ref{tab:NAS-CIFAR10}, the best values achieved by PWAS and PWASp after 100 evaluations are often already better or comparable to the results obtained by some other algorithms after 200 evaluations. These observations are also reflected in
Figures~\ref{fig:convergence_plot}--\ref{fig:data_profile}, where we observe that PWAS, PWASp, and MISO consistently demonstrate superior early-stage performance. For $\tau = 0.5$, they can solve all the benchmark problems with fewer than 10 equivalent simplex gradients. This suggests that these three algorithms are well-suited for applications where rapid, early-stage reduction in function value is critical. However, for $\tau = 0.1$, the performance of PWASp quickly drops.
It is also observed that PWAS performs consistently better than PWASp. It is because that PWAS has access to function evaluations, while PWASp only receives pairwise comparisons, making it struggle to obtain highly accurate results. Nonetheless, in spite of the more limited information it gets, PWASp outperforms several other solvers (One-hot BO, SMAC, TPE, and EXP3BO) in most of the tested benchmarks. 
We also stress that PWASp is advantageous when the objective function is not easily quantifiable, but pairwise comparisons can be made. 
On the other hand, When $\tau$ decreases, the performances of PWAS and MISO remain in the top three, with the performance of MISO drops slightly.
CoCaBO-0.5 also holds a competitive position but is not as dominant as PWAS in early performance ratio (up to 2). However, as the performance ratio increases, the performance of CoCaBO-0.5 remains steady. Additionally, CoCaBO-0.5 showed superior performance in both real-world problems (\textsf{XG-MNIST} and \textsf{NAS-CIFAR10}). Although, it is worth noting that the performance differences are rather small in scale.

\begin{figure}[tbh!]
    \centering
    \includegraphics[width=.95\textwidth]{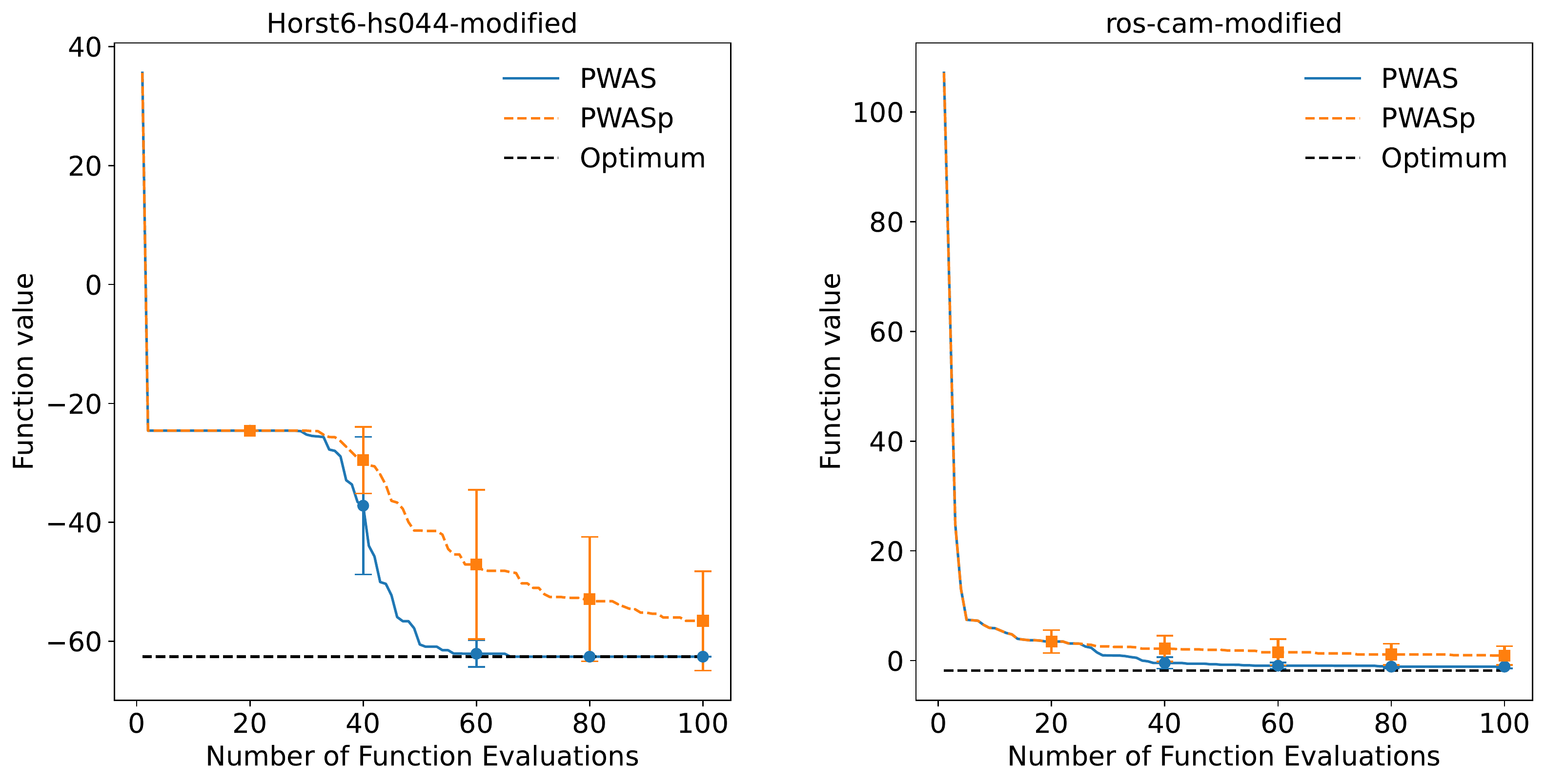}\quad\quad 

    \caption{Convergence graphs for the constrained benchmark problems. Error bars indicate standard deviations, which are obtained with 20 random repetitions.}
    \label{fig:convergence_constrained}
\end{figure}

The optimization results for the constrained benchmarks 
are shown in Table~\ref{tab:constrained_synthetic} with the convergence graphs shown in Figure~\ref{fig:convergence_constrained}. We observe that both PWAS and PWASp can quickly reduce the objective function values after a small number of black-box function evaluations, demonstrating their ability to handle mixed-variable linear QUAK constraints Also, PWAS achieves better results than PWASp regarding the best values obtained after the maximum allowed black-box function evaluations and consistency over multiple repetitions. 

\begin{table}[htb!]
\caption{Performance of PWAS/PWASp on constrained mixed-variable synthetic problems.}
\centering
\begin{tabularx}{\linewidth}{>{\hsize=1.4\hsize}X *{4}{>{\hsize=0.9\hsize}X}}
\hline
\multirow{2}{*}{Algorithm}  & \multicolumn{2}{l}{\textsf{Horst6-hs044-modified}}   & \multicolumn{2}{l}{\textsf{ros-cam-modified}}  \\ \cline{2-5} 
                                     & mean           & std           & mean           & std          \\ \hline
\textsf{PWAS}           & -62.579                  & 3.5275e-08              & -1.1151                & 0.3167               \\
\textsf{PWASp}          & -56.5539                 & 8.3454                  & 0.90540                 & 1.7200             \\ \hline
Global optimum & \multicolumn{2}{l}{-62.579}                        & \multicolumn{2}{l}{-1.81}                     \\ \hline
\multicolumn{5}{l}{\fontsize{8}{9}\selectfont Note: for both PWAS and PWASp, the optimum is obtained after 100 black-box    }\\
\multicolumn{5}{l}{\fontsize{8}{9}\selectfont function evaluations. Statistics are obtained with 20 random  repetitions.}\\
\end{tabularx}
\label{tab:constrained_synthetic}
\end{table}

To show the efficiency of PWAS and PWASp, the average CPU time spent by PWAS and PWASp to fit the surrogate and solve the acquisition problem to suggest the next sample to query is reported in Table~\ref{tab:cpu_time}
for each tested benchmark problem. Considering that often evaluating the black-box function $f$ or comparing samples involves expensive-to-evaluate simulations or experiments, such a CPU time can be considered negligible in real-life applications.

\begin{table}[h!]
\caption{CPU time (s) for surrogate fitting and acquisition optimization, averaged over $N_{\rm max}- N_{\rm init}$ active sampling steps.}
\centering
\begin{tabularx}{\textwidth}{>{\hsize=1.4\hsize}X *{8}{>{\hsize=0.95\hsize}X}}
\hline
 & \multicolumn{1}{c}{}                                    & \textsf{Func-2C} & \textsf{Func-3C} & \textsf{Ackley-5C} & \textsf{XG-MNIST} & \textsf{NAS-CIFAR10} & \textsf{Horst6-hs044-modified} & \textsf{ros-cam-modified} \\ \hline
\multirow{2}{*}{\begin{tabular}[c]{@{}l@{}}Surrogate\\ fitting\end{tabular}}                 & \textsf{PWAS}                                           & 0.565            & 0.556            & 0.889              & 0.327             & 0.329                & 0.211                          & 0.198                     \\
& \textsf{PWASp}                                          & 0.221            & 0.289            & 0.544              & 0.312             & 0.422                & 0.177                          & 0.162                     \\ \hline
\multirow{2}{*}{\begin{tabular}[c]{@{}l@{}}Acquisition\\ optimization\end{tabular}}          & \textsf{PWAS}                                           & 0.231            & 0.196            & 1.250              & 0.505             & 1.871                & 0.327                          & 0.311                     \\
& \textsf{PWASp}                                          & 0.270            & 0.420            & 1.352              & 0.589             & 1.700                & 0.387                          & 0.364                     \\ \hline                                        \end{tabularx}
\label{tab:cpu_time}
\end{table}

\section{Conclusion}\label{sec:conclusion}
The algorithms PWAS and PWASp introduced in this paper can handle global and preference-based optimization problems involving mixed variables subject to linear preference-based optimization problems involving mixed variables subject to linear quantifiable unrelaxable a priori known (QUAK) constraints. Tests on different synthetic and real-world benchmark problems show that PWAS and PWASp can obtain better or comparable performance than other existing methods. 
In addition, the proposed acquisition strategies in PWAS and PWASp does not require uncertainty measurements, which can be extended to other simpler surrogate models such as polynomials.
Although convergence to global optimizers cannot be guaranteed, we observed that PWAS and PWASp could quickly reduce (increase for maximization problems) the objective function values within a limited number of black-box function evaluations, despite the presence of integer and categorical variables and mixed-integer linear constraints. Therefore, both PWAS and PWASp can be considered as good heuristic algorithms for mixed-variable black-box optimization problems.

Future research will be devoted to extend this approach to handle mixed-variable problems under mixed-variable nonlinear constraints. Several approaches can be used, such as replacing them by piecewise affine approximations. Additionally, it can be interesting to integrate the proposed exploration function with other surrogate models (e.g., tree-based BO and polynominals).

\newpage
\noindent \textbf{Statements and Declarations}\\

\noindent\textbf{\textit{Funding}}\\

\noindent The authors declare that no funds, grants, or other support were received during the preparation of this manuscript.

\vspace{3em}
\noindent\textbf{\textit{Competing Interests}}\\

\noindent The authors have no relevant financial or non-financial interests to disclose 
\vspace{3em}

\noindent\textbf{\textit{Data Availability}}\\

\noindent The PWAS package and the tested benchmarks are available on the GitHub repository (\url{https://GitHub.com/mjzhu-p/PWAS}). The MNIST dataset analysed for benchmark \textsf{XG-MNIST} in this study is publicly available and retrieved via \textsf{sklearn.datasets.load\_digits}. The NASBench dataset analysed for benchmark \textsf{NAS-CIFAR10} is publicly available from the repository hosted on GitHub (\url{https://GitHub.com/google-research/nasbench}) under Apache License 2.0 and can be downloaded via \url{https://storage.googleapis.com/nasbench/nasbench_only108.tfrecord}. The MISO package is publicly available from the repository  hosted on GitHub (\url{https://github.com/Julie2901/miso}). The NOMAD4 package is publicly available from the repository hosted on GitHub  (\url{https://github.com/bbopt/nomad}) under LGPL-3.0 license.

\newpage
\begin{appendices}

\section{Illustrative Example - Surrogate Fitting}~\label{sec:app_ill3}
\begin{figure}[tbh]
    \centering
    \includegraphics[width=\textwidth]{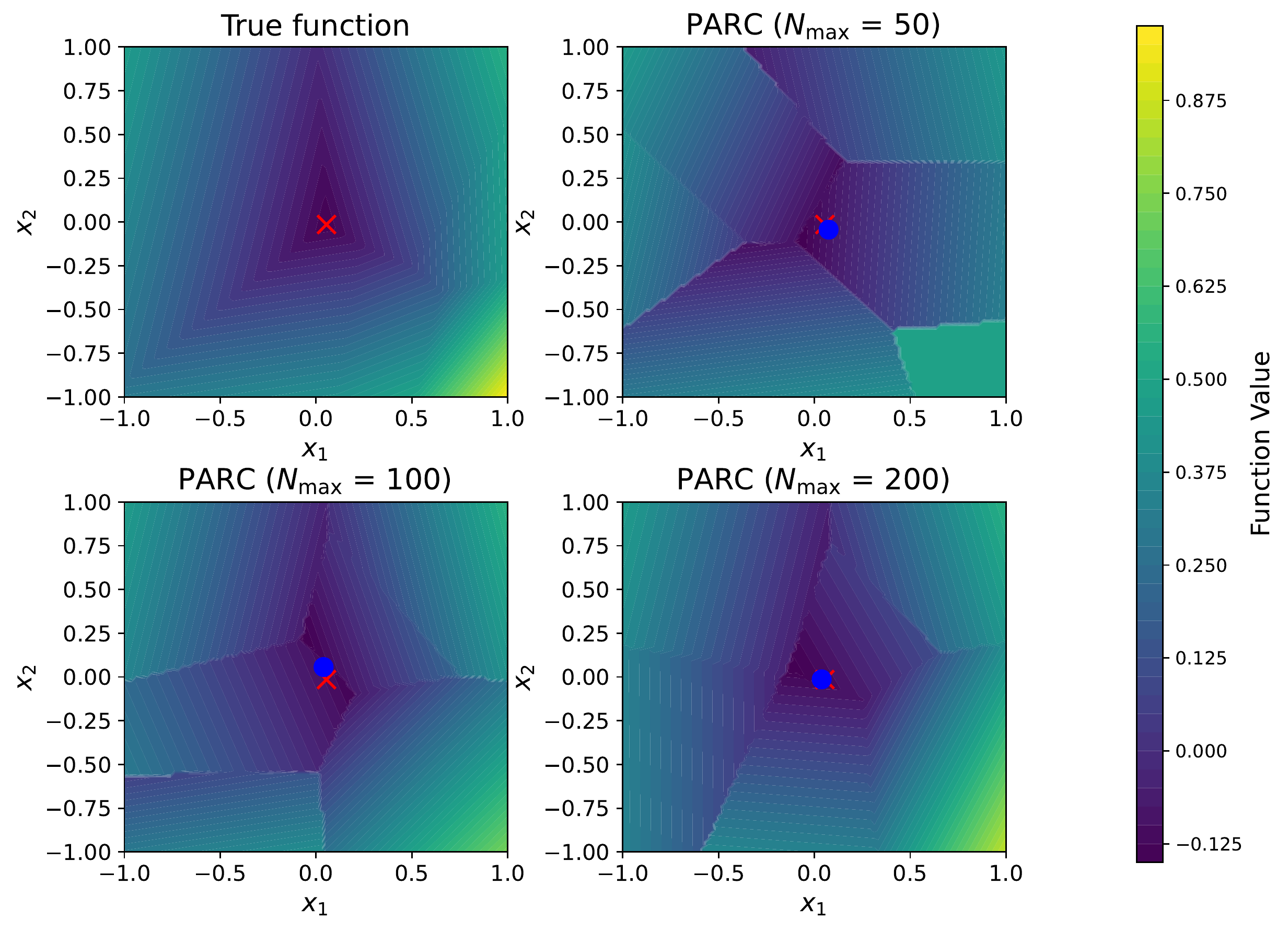}
    \caption{Contour plots of Function~\eqref{eq:PWA-example-fun}, with the ground truth (top-left) and approximations generated by PWA surrogates with different numbers of samples ($N_{\rm max}$). The approximations are based on 50, 100, and 200 samples. Red cross: optimum; blue dot: best sample obtained within $N_{\rm max}$ function evaluations. }
    \label{fig:contour_ill3}
\end{figure}

\begin{figure}[tbh]
    \centering
    \includegraphics[width=\textwidth]{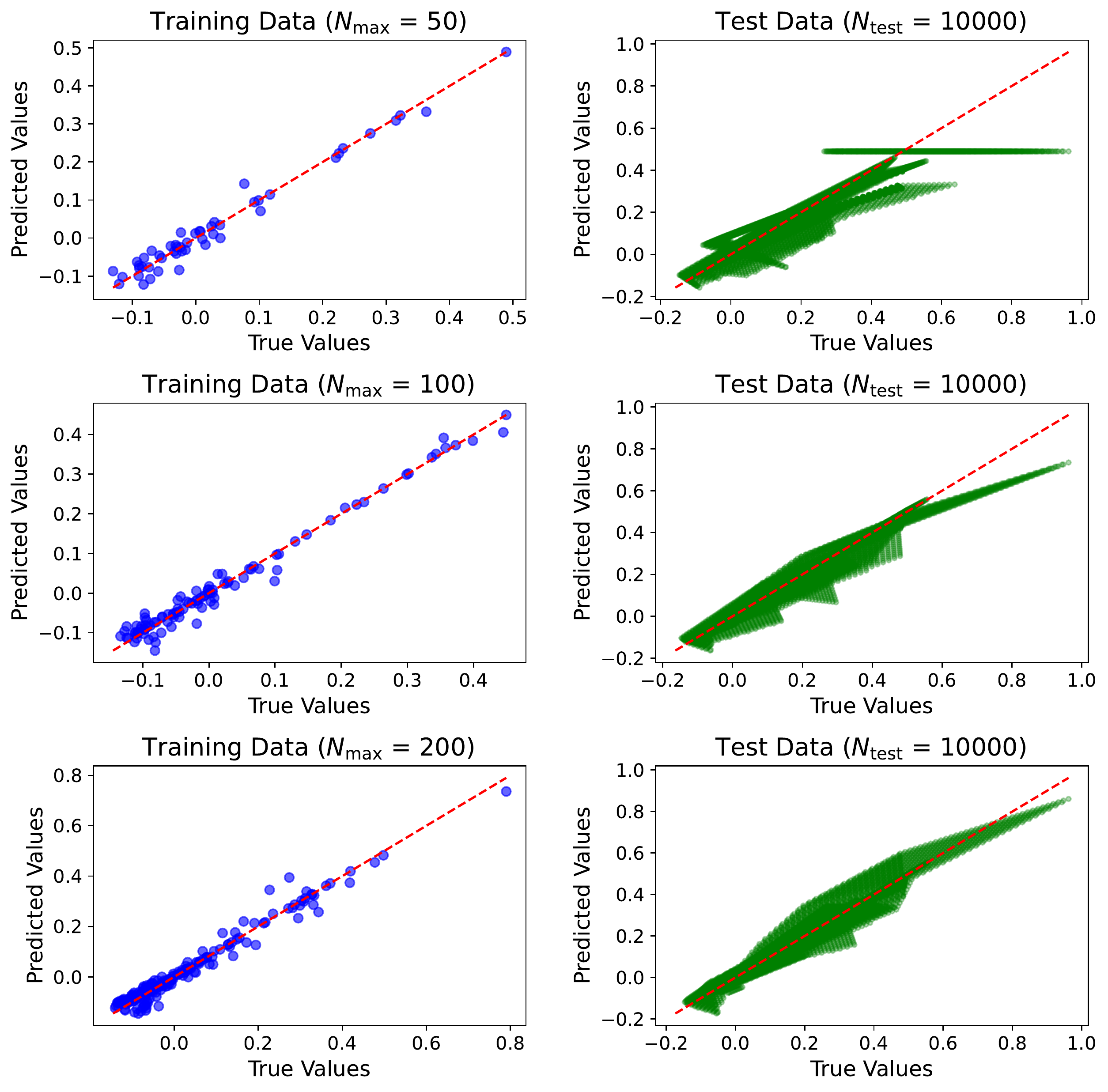}
    \caption{Comparisons between the true function values of~\eqref{eq:PWA-example-fun} and the predicted function values using the PWA surrogates for both training (with different values of $N_{\rm max}$) and testing samples (with 10,000 samples). Blue dots: training samples collected during optimization; green dots: testing samples; red dashed lines: ideal case where true values equal to the predicted values.}
    \label{fig:tru_vs_predict_ill_3}
\end{figure}

To illustrate the remark noted at the end of Section~\ref{subsec:PWA_surrogate}, we optimize the following PWA function from~\cite{BemPARC21},
\begin{eqnarray}
    f(x)&=&\max\left\{
    \smallmat{0.8031\\0.0219\\-0.3227}'\smallmat{x_1\\x_2\\1}
    ,
    \smallmat{0.2458\\-0.5823\\-0.1997}'\smallmat{x_1\\x_2\\1}
    ,
    \smallmat{0.0942\\-0.5617\\-0.1622}'\smallmat{x_1\\x_2\\1}
    ,
    \smallmat{0.9462\\-0.7299\\-0.7141}'\smallmat{x_1\\x_2\\1}
    ,
    \right.\nonumber\\&&\left.
    \smallmat{-0.4799\\0.1084\\-0.1210}'\smallmat{x_1\\x_2\\1}
    ,
    \smallmat{0.5770\\0.1574\\-0.1788}'\smallmat{x_1\\x_2\\1}
    \right\}.
    \label{eq:PWA-example-fun}
\end{eqnarray}
We employ PWAS with varying maximum numbers of function evaluations $N_{\rm max}$, using the samples collected during optimization to fit~\eqref{eq:PWA-example-fun}. As shown in Figure~\ref{fig:contour_ill3}, after 50 function evaluations, the surrogate effectively identifies the optimal region (dark purple region in Figure~\ref{fig:contour_ill3}), and after 200 evaluations, it closely approximates the optimal region. 
Furthermore, Figure~\ref{fig:tru_vs_predict_ill_3} demonstrates that as $N_{\rm max}$ increases from 50 to 200, the surrogate's predictions align more closely with the true values across the test dataset. With a smaller $N_{\rm max}$, predictions are most accurate near the optimum; as $N_{\rm max}$ grows, the surrogate's accuracy improves also in regions further from the optimum. This indicates that PARC can be effectively used for our purpose as noted in the remark.

\noindent Note: default solver settings of PWAS are used when collecting fitting samples for the illustrative example.

\newpage

\section{Benchmark}\label{sec:app_bench}
\noindent Note: the unconstrained mixed-variable synthetic and real-world problems are adopted from~\cite{ru2020bayesian} for our comparisons. Note that the objective function is maximized in~\cite{ru2020bayesian}, so we consider the minimization of $-f(X)$ in PWAS and PWASp.
\subsection{Unconstrained Mixed-variable Synthetic Benchmarks}~\label{subsec:uncont_synthetic}

\noindent \textsf{Func-2C}~\cite{ru2020bayesian,dixon1978global,simulationlib,molga2005test,beale1958iterative,more1981testing}: 
$n_c = 2$, $n_{\rm int} = 0$, and $n_d = 2$ with $n_{i} = 3$ for each categorical variable (denoted as $n_{di}$, for $i = 1, 2$). Each categorical variable is in $\{0, 1, 2\}$. The bounds are $\ell_x = [-1.0\ -1.0]\tr$, $u_x = [1.0\ 1.0]\tr$. The global maximum $f(X) = 0.20632$ is
attained at $X = [0.0898\ -0.7126\ 1\ 1]\tr$ and $[-0.0898\ 0.7126\ 1\ 1]\tr$.\\
\begin{equation}\label{eq:Func-2C}
\begin{split}
	f(X) &=  \begin{cases}
 				f_1 + f_{\rm ros}(x)  &  n_{d2}= 0\\
				f_1 + f_{\rm cam}(x)  &n_{d2}= 1 \\
				f_1 + f_{\rm bea}(x)  &n_{d2}= 2, 
			\end{cases} \\
	{\rm where\ }  
	&\quad f_1(x,y) =\begin{cases}
					f_{\rm ros}(x) & n_{d1}= 0\\
					f_{\rm cam}(x) & n_{d1}= 1 \\
					f_{\rm bea}(x) & n_{d1}= 2 
			  \end{cases} \\
	&\quad f_{\rm ros}(x) = -(100 (x_2 - x_1^2)^2 + (x_1 - 1)^2)/300\\
	 &\quad f_{\rm cam}(x) = -(a_1 + a_2 + a_3)/10\\
	 	&\quad \quad \quad \quad a_1 = (4 - 2.1 x_1^2 + \frac{x_1^4}{3}) x_1^2 \\
	 	&\quad \quad \quad \quad a_2 = x_1 x_2 \\
	&\quad \quad \quad \quad a_3 = (-4 + 4 x_2^2) x_2^2 \\
	&\quad f_{\rm bea}(x) = -((1.5 - x_1 + x_1 x_2)^2 + (2.25 - x_1 + x_1 x_2^2)^2 \\
	&\quad \quad \quad \quad \quad \quad + (2.625 - x_1 + x_1 x_2^3)^2)/50\\
\end{split}
\end{equation}


\newpage
\noindent \textsf{Func-3C}~\cite{ru2020bayesian,dixon1978global,simulationlib,molga2005test,beale1958iterative,more1981testing}: 
$n_c = 2$, $n_{\rm int} = 0$, and $n_d = 3$ with $n_{i} = 3$ for each categorical variable (denoted as $n_{di}$, for $i = 1, 2, 3$).
Each categorical variable is in $\{0, 1, 2\}$. The bounds are $\ell_x = [-1.0\ -1.0]\tr$, $u_x = [1.0\ 1.0]\tr$. The global maximum
$f(X) = 0.72214$ is attained at $X = [0.0898\ -0.7126\ 1\ 1\ 0]\tr$ and $[-0.0898\ 0.7126\ 1\ 1\ 0]\tr$. \\
\begin{equation}\label{eq:Func-3C}
\begin{split}
	f(X) &=  \begin{cases}
 				f_2 + 5f_{\rm cam}(x)  &  n_{d3}= 0\\
				f_2 + 2f_{\rm ros}(x)  &n_{d3}= 1 \\
				f_2 + n_{d2}f_{\rm bea}(x)  &n_{d3}= 2, 
			\end{cases} \\
	{\rm where\ } 
	&\quad f_2(x,y) =\begin{cases}
					f_1 + f_{\rm ros}(x) & n_{d21}= 0\\
					f_1 + f_{\rm cam}(x) & n_{d2}= 1 \\
					f_1 + f_{\rm bea}(x) & n_{d2}= 2 
			  \end{cases} \\ 
	&\quad f_1(x,y) =\begin{cases}
					f_{\rm ros}(x) & n_{d1}= 0\\
					f_{\rm cam}(x) & n_{d1}= 1 \\
					f_{\rm bea}(x) & n_{d1}= 2 
			  \end{cases} \\
	&\quad f_{\rm ros}(x), f_{\rm cam}(x), {\rm \ and\ } f_{\rm bea}(x) {\rm \ are\  defined\  in}~\eqref{eq:Func-2C}
\end{split}
\end{equation}

\vspace{2em}

\noindent \textsf{Ackley-5C}~\cite{ru2020bayesian,simulationlib,ackley1987model}:
$n_c = 1$, $n_{\rm int} = 0$, and $n_d = 5$ with $n_{i} = 17$ for each category (denoted as $n_{di}$, for $i = 1, \dots, 5$). Each categorical variableis in $\{0, 1, \dots, 16\}$. The bounds are $\ell_x = -1.0 $, $u_x = 1.0$. The global maximum $f(X) = 0$
is attained at $X = [0\  8\ 8\ 8\ 8\ 8]\tr$. \\
\begin{equation}\label{eq:Ackley-5C}
\begin{split}
	f(X) &= a \exp\left(-b\left(\frac{s_1}{n}\right)^{\frac{1}{2}}\right) + \exp\left(\frac{s_2}{n}\right) - a - \exp(1),\\
	{\rm where\ } 
	&\quad  a = 20, b = 0.2, c = 2\pi \\
	&\quad  s_1 = x^2 + \sum_{i = 1}^{5} z_i^2(n_{di})\\
	&\quad  s_2 = \cos(cx) + \sum_{i = 1}^{5} \cos (c z_i(n_{di}))\\
	&\quad  z_i(n_{di}) = -1 + 0.125n_{di}
\end{split}
\end{equation}

\vspace{1em}
\subsection{Unconstrained Mixed-variable Real-world Benchmarks}~\label{subsec:uncont_real}

\noindent\textsf{XG-MNIST}~\cite{ru2020bayesian,chen2016xgboost,lecun2010mnist}:
$n_c = 4$, $n_{\rm int} = 1$, and $n_d = 3$ with $n_{i} = 2$, for $i = 1, 2, 3$. Each categorical variable ($n_{di}$) can be either 0 or 1. The bounds are $\ell_x = [10^{-6},10^{-6}\ 0.001\ 10^{-6}]\tr$, $u_x = [1\ 10\ 1\ 5]\tr$; $\ell_y = 1$, $u_y = 10$.\\

\noindent \underline{Notes on the optimization variables}:\\
The 0.7/0.3 stratified train/test split ratio is applied as noted in~\cite{ru2020bayesian}. The \emph{xgboost} package is used~\cite{chen2016xgboost} on MNIST classification~\cite{lecun2010mnist}. 
The optimization variables in this problem are the parameters of the \emph{xgboost} algorithm.
Specifically, the continuous variables $x_1$, $x_2$, $x_3$, and $x_4$ refer to the following parameters in \emph{xgboost}, respectively: \texttt{`learning\_rate'}, \texttt{`min\_split\_loss'}, \texttt{`subsample'}, and \texttt{`reg\_lambda'}. The integer variable $y$ stands for the \texttt{`max\_depth'}. As for the categorical variables, $n_{d1}$ indicates the booster type in \emph{xgboost}, where $n_{d1} = \{0, 1\}$ corresponding to \{\texttt{`gbtree'}, \texttt{`dart'}\}. $n_{d2}$ represents the \texttt{`grow\_policy'}, where $n_{d2} = \{0, 1\}$ corresponding to \{\texttt{`depthwise'}, \texttt{`lossguide'}\}. $n_{d3}$ refers to the \texttt{`objective'}, where $n_{d3} = \{0, 1\}$ corresponding to \{\texttt{`multi:softmax'}, \texttt{`multi:softprob'}\}. \\

\noindent \underline{Notes on the objective function}:\\
The classification accuracy on test data is used as the objective function.

\vspace{3em}
\noindent\textsf{NAS-CIFAR10}~\cite{ru2020bayesian,ying2019bench,klein2019tabular}:
$n_c = 21$, $n_{\rm int} = 1$, and $n_d = 5$ with $n_{i} = 3$ for each category (denoted as $n_{di}$ for $i = 1, \dots, 5$). Each categorical variable is in $\{0, 1, 2\}$. The bounds are $\ell_x^i = 0$, $u_x^i = 1$, $\forall i=1,\ldots,{n_c}$; 
$\ell_y = 0$, $u_y = 9$.\\

\noindent \underline{Notes on the optimization problem}:\\
The public dataset NAS-Bench-101~\cite{yang2019machine} is used. This dataset maps convolutional neural network (CNN) architectures to their trained and evaluated performance on CIFAR-10 classification. As a result, we can quickly look up the validation accuracy of the proposed CNN architecture by PWAS/PWASp. The same encoding method for the CNN architecture topology as noted in~\cite{yang2019machine,ru2020bayesian} is used. \\

\noindent \underline{Notes on the optimization variables}:\\
The CNN architecture search space is described by a directed acyclic graph (DAG) which has 7 nodes with the first and the last nodes being the input and output nodes. 
The continuous variables represent the probability values for the 21 possible edges in the DAG. The integer variable $y$ is the number of edges present in the DAG. The categorical variables are the operations for the 5 intermediate nodes in the DAG, for which $n_{di} = \{0, 1, 2\}$ corresponding to \{\texttt{`3x3 conv'}, \texttt{`1x1 conv'}, \texttt{`3x3 max-pool'}\}. Within the 21 possible edges, only $y$ edges with the highest probability are activated. The DAG that leads to invalid CNN architecture topology specifications will result in zero validation accuracy.\\

\noindent \underline{Notes on the objective function}:\\
The validation accuracy of the defined CNN architecture topology on CIFAR-10 classification is used as the objective function.

\newpage
\subsection{Constrained Mixed-variable Synthetic Problems}~\label{subsec:conts_synthetic}

\noindent\textsf{Horst6-hs044-modified}~\cite{horst2000introduction,le2012optimizing}:
$n_c = 3$, $n_{\rm int} = 4$, and $n_d = 2$ with $n_{1} = 3$ and $n_{2} = 2$ (the first ($n_{d1}$) and the second ($n_{d2}$) categorical variable are in $\{0, 1, 2\}$ and $\{0, 1\}$, respectively. $\ell_x = [0\ 0\ 0]\tr$, $u_x = [6\ 6\ 3]\tr$; $\ell_y = [0\ 0\ 0\ 0]\tr$, $u_y = [3\ 10\ 3\ 10]\tr$. The global minimum $f(X) =-62.579$ is attained at $X = [5.21066\ 5.0279\ 0\ 0\ 3\ 0\ 4\ 2\ 1]\tr$.\\
\begin{equation}\label{eq:Horst6-hs044-modified}
\begin{split}
	f(X) &=  \begin{cases}
 				| f_1 |   &  n_{d2}= 0\\
				f_1 &n_{d2}= 1
			\end{cases} \\
	{\rm s.t\ }  
	&\quad A_{\rm ineq}x + B_{\rm ineq}y \leq b_{\rm ineq},\\	
	{\rm where\ }  
	&\quad f_1(x,y) =\begin{cases}
					f_{\rm Horst6}(x) +  f_{\rm hs044}(y) & n_{d1}= 0\\
					0.5f_{\rm Horst6}(x) +  f_{\rm hs044}(y) & n_{d1}= 1\\
					f_{\rm Horst6}(x) +  2f_{\rm hs044}(y) & n_{d1}= 2\\
			  \end{cases} \\
	&\quad f_{\rm Horst6}(x) = x^TQx + px\\
	&\quad \quad \quad \quad Q = \begin{bmatrix}0.992934& - 0.640117& 0.337286 \\-0.640117& - 0.814622& 0.960807\\ 0.337286& 0.960807 & 0.500874\end{bmatrix} \\
	&\quad \quad \quad \quad p =  \begin{bmatrix}-0.992372& - 0.046466& 0.891766 
	 \end{bmatrix} \\
	 &\quad f_{\rm hs044}(y) = x_0 - x_1 - x_2 - x_0x_2 + x_0x_3 + x_1x_2 - x_1x_3\\
	&\quad  A_{\rm ineq} =  \begin{bmatrix}
0.488509& 0.063565& 0.945686 \\
-0.578592& -0.324014& -0.501754 \\
-0.719203& 0.099562& 0.445225\\
-0.346896& 0.637939& -0.257623\\
-0.202821& 0.647361& 0.920135\\
-0.983091& -0.886420& -0.802444\\
-0.305441& -0.180123& -0.515399
\end{bmatrix} \\
	 &\quad B_{\rm ineq} =  \begin{bmatrix}
1& 2& 0& 0 \\
4& 1& 0& 0\\
3& 4& 0& 0\\
0& 0& 2& 1\\
0& 0& 1& 2\\
0& 0& 1& 1
\end{bmatrix}  \\
	 &\quad b_{\rm ineq} = [2.86506,\  -1.49161,\   0.51959,\   1.58409,\   2.19804,\   -1.30185,\\
	 &\quad\quad \quad \quad\   -0.73829,\   8,\   12,\   12,\   8,\   8,\   5]\tr
\end{split}
\end{equation}

\newpage
\noindent\textsf{ros-cam-modified}~\cite{jamil2013literature,Bem20}:
$n_c = 2$, $n_{\rm int} = 1$, and $n_d = 2$ with $n_{i} = 2$ for each categorical variable (denoted as $n_{di}$, for $i = 1, 2$). Each categorical variable is in $\{0, 1\}$. $\ell_x = [-2.0\ -2.0]\tr$, $u_x = [2.0\ 2.0]\tr$; $\ell_y = 1$, $u_y = 10$. The global minimum $f(X) =-1.81$ is attained at $X = [0.0781\ 0.6562\ 5\ 1\ 1]\tr$.\\
\begin{equation}\label{eq:ros-cam-modified}
\begin{split}
	f(X) &=  \begin{cases}
 				f_1 + f_{\rm ros}(x,y)  &  n_{d2}= 0\\
				f_1 + f_{\rm cam}(x,y)  &n_{d2}= 1
			\end{cases} \\
	{\rm s.t\ }  
	&\quad A_{\rm ineq}x \leq b_{\rm ineq},\\	
	{\rm where\ }  
	&\quad f_1(x,y) =\begin{cases}
					f_{\rm ros}(x,y) & n_{d1}= 0\\
					f_{\rm cam}(x,y) & n_{d1}= 1
			  \end{cases} \\
	&\quad f_{\rm ros}(x,y) = 100 (x_2 - x_1^2)^2 + (x_1 - 1)^2 + (y-3)^2\\
	 &\quad f_{\rm cam}(x,y) = a_1 + a_2 + a_3 + (y-5)^2\\
	 	&\quad \quad \quad \quad a_1 = (4 - 2.1 x_1^2 + \frac{x_1^4}{3}) x_1^2 \\
	 	&\quad \quad \quad \quad a_2 = x_1 x_2 \\
	&\quad \quad \quad \quad a_3 = (-4 + 4 x_2^2) x_2^2 \\
	&\quad  A_{\rm ineq} =  \begin{bmatrix}
1.6295& 1 \\
0.5&3.875 \\
-4.3023& -4\\
-2&1\\
0.5&-1
\end{bmatrix} \\
	 &\quad b_{\rm ineq} = [3.0786,\  3.324,\  -1.4909,\  0.5,\   0.5]\tr
\end{split}
\end{equation}

\newpage
\section{Algorithms}\label{sec:app_algorithm}
\subsection{CoCaBO}
Continuous and Categorical Bayesian Optimization (CoCaBO)~\cite{ru2020bayesian} is an algorithm proposed to optimize box-constrained expensive black-box problems with both continuous and categorical variables, specifically for problems with \emph{multiple} categorical variables with \emph{multiple} possible values. CoCaBO model the input space with a Gaussian Process kernel, which is designed to allow information sharing across different categorical variables to enhance data efficiency. 

\subsection{EXP3BO}
EXP3BO~\cite{gopakumar2018algorithmic} is an algorithm proposed to optimize box-constrained, expensive black-box problems, which is modified from BO via the EXP3 algorithm~\cite{auer2002nonstochastic,seldin2013evaluation}. It can deal with mixed categorical and continuous input spaces by utilizing multi-armed bandits. Specifically, EXP3BO constructs a Gaussian Process surrogate specific to each chosen category, making it unsuitable for handling problems with multiple categorical classes.

\subsection{One-hot BO}
One-hot BO~\cite{gpyopt2016} handles the input space with continuous and categorical variables by one-hot encoding the categorical variables and treating these encoded variables as continuous. Then, the standard BO procedure is used to solve the optimization problem on the transformed input space.

\subsection{TPE}
Tree-structured Parzen Estimator (TPE)~\cite{bergstra2011algorithms} is a black-box optimization algorithm based on tree-structured Parzen density estimators. TPE uses nonparametric Parzen kernel density estimators to model the distribution of good and bad configurations w.r.t. a reference value. Due to the nature of kernel density estimators, TPE also supports continuous and categorical spaces.

\subsection{MISO}
Mixed-Integer Surrogate Optimization (MISO)~\cite{muller2016miso} is an algorithm that targets expensive black-box functions with mixed-integer variables. It constructs a surrogate model using the radial basis function to approximate the unknown function. MISO follows a general procedure of surrogate-based optimization methods. Additionally, MISO combines different sampling strategies (\emph{e.g.}, coordinate perturbation, random sampling, expected improvement, target value, and surface minimum) and local search to obtain high-accuracy solutions.

\subsection{NOMAD}
NOMAD~\cite{abramson2011nomad,audet2021nomad} is a C++ implementation of the Mesh Adaptive
Direct search (MADS). It is designed to solve difficult black-box optimization problems. In particular, it can handle nonsmooth, nonlinearly constrained, single or bi-objetive, and mixed variable optimization problems. It handles the categorical variable using the extended poll, which is defined as following~\cite{le2019nomad}:
\begin{quote}
    The extended poll first calls the user-provided
procedure defining the neighborhood of categorical variables. The procedure returns a list of points that are neighbors of the current best point (incumbent) such that categorical variables
are changed and the other variables may or may not be changed. These points are called the
extended poll points and their dimension may be different than the current best point, for
example when a categorical variable indicates the number of continuous variables.
\end{quote}

\subsection{SMAC}
Sequential Model-based Algorithm Configuration (SMAC)~\cite{hutter2011sequential} is a surrogate-based black-box optimization method originally proposed to tackle \emph{algorithm configuration} problems with continuous and categorical variables. Its model is based on random forests~\cite{breiman2001random}, so it can handle categorical variables explicitly. SMAC uses empirical mean and variance within a tree ensemble to identify uncertain search space regions and optimizes the acquisition function by combing local and random search. 

\end{appendices}

\newpage
\bibliographystyle{spmpsci} 
\bibliography{Biblio_pwa_pref}  


\end{document}